\documentstyle[11pt]{article}

\def\proof{\noindent{\bf Proof}$\;\;$}      
\def\define{\bigskip\noindent{\bf Definition}$\;\;$} 
      
\def\hook{\,\hbox to 10pt{\vbox{\vskip 6pt\hrule width 6.5pt height 1pt}
         \kern -4.0pt\vrule height 8pt width 1pt\hfil}\,}
\def\blob{\mbox{$\;\Box$}}
\def\beq{\begin{equation}}
\def\eeq{\end{equation}}
\def\beqa{\begin{eqnarray}}
\def\eeqa{\end{eqnarray}}
\def\beqann{\begin{eqnarray*}}
\def\eeqann{\end{eqnarray*}}
\def\vf{{{VF}(Y)}}
\def\tauvf{\tau^*\vf}
\def\pttauvf{(\pi_{XY}(y),y,f_i,\epsilon_A)}
\def\lvy{L_VY}
\def\ptlvy{{(y,\{e_i,\epsilon_A\})}}

\def\frame{\{e_i,\epsilon_A\}}
\def\phibl{{\phi_{(B,\lambda)}}}
\def\projectable{{\cal X}_{\mbox{\tiny Proj}}Y}
\def\wtheta #1{{\wedge^{#1}\theta}}
\def\iwtheta #1{{\wedge^{#1}i^*\theta}}          
\def\tvlvy{T^1_V(\lvy)}
\def\xhat #1{{X_{\hat {#1}}}}
\def\vectorclass #1{{\lbrack\!\lbrack X_{\hat #1}\rbrack\!\rbrack}}
\def\pb#1#2{\{\hat #1,\hat #2\}}
\def\r{\mbox{$\bf R$}}

\def\rnplusk{\r^{n+k}}
\def\rnk{{\r^{n\times k}}}
\def\rkn{{\r^{k\times n}}}
\def\rnkxr{{\r^{n\times k}\times \r}}
\def\nka{ \left( \begin{array}{cc} N&0\\ A&K \end{array}\right)}

\def\f{{\boldmath \mbox{$f$}}}
\def\lin{\mbox{Lin$\,$}}
\renewcommand{\hom}{\mbox{Hom$\,$}}
\renewcommand{\Pr}{\mbox{Pr$\,$}}
\def\basisx #1{{\partial\over\partial x^{#1}}}
\def\basisy #1{{\partial\over\partial y^{#1}}}
\def\basispi #1 #2{{\partial\over\partial \pi^{#1}_{#2}}}
\def\basisp #1 #2{{\partial\over\partial p^{#1}_{#2}}}
\def\partialx #1 #2{{\partial #1 \over \partial x^{#2}}}
\def\partialy #1 #2{{\partial #1 \over \partial y^{#2}}}
\def\xhatf{{X_{\hat f}}}
\def\multii#1{i_1\dots i_{#1}}
\def\ga{{G_{\cal A}}}
\def\lb2{[\![}
\def\rb2{]\!]}
\def\tr{\mbox{Tr}\,}
\def\gunther{\hom_Y(V(TY),TX)}
\def\kt{\hom_Y(V(TY),\wedge^{n-1}X)}

\newtheorem{thm}{Theorem} [section]
\newtheorem{lemma}[thm]{Lemma}            
\newtheorem{cor}[thm]{Corollary}         
\newcounter{ize}

\small\normalsize

\title{A Frame Bundle Generalization
of Multisymplectic Field Theories}

\author{J. K. Lawson\thanks{Electronic mail: lawson@mthcsc.wfu.edu}\\
Department  of Mathematics and Computer Science\\
 Wake Forest University\\
Winston-Salem, NC 27109-7388}
\date{31 May 1997}

\begin{document}

\maketitle

\begin{abstract} 
This paper presents a generalization of symplectic 
 geometry 
to a principal bundle over
 the configuration space of a classical field.
This bundle, the vertically adapted linear frame bundle, 
is obtained by breaking the symmetry
 of the full linear frame bundle of the field configuration space, and it 
inherits a generalized symplectic structure from the full frame bundle.
The geometric structure of the vertically adapted frame bundle admits
vector-valued field observables and 
produces vector-valued
Hamiltonian vector fields, from which we can define a Poisson bracket
on the field observables. 
We show that the linear and affine multivelocity spaces and multiphase spaces 
for geometric field theories
are associated to the vertically adapted frame bundle.
In addition, the new geometry not only generalizes 
both the linear and the affine models of multisymplectic geometry
but also resolves fundamental
problems found in both mul\-ti\-sym\-plec\-tic models.
\end{abstract}
\thispagestyle{empty}

\bigskip

\noindent{\em Keywords: symplectic geometry, mul\-ti\-sym\-plec\-tic 
geometry, frame bundle, Hamiltonian field theories, Poisson bracket.

\medskip

\noindent MS classification:  58 F 05 (Primary)  70 G 50 (Secondary)

\medskip

 }

\section{\normalsize\bf  INTRODUCTION}
\setcounter{equation}{0}

Norris's theory of $n$-symplectic geometry${}^{\ref{No1},\ref{No2}}$ of the linear frame bundle
has proven to be a powerful tool in
the study of symplectic geometry and  and the standard geometric
model of Hamiltonian par\-ti\-cle mechanics.
In Norris's theory, the canonical
soldering form of the linear frame bundle of a manifold of particle
configurations behaves as a
vector-valued {\em $n$-symplectic}
potential. 
The ensuing
$n$-symplectic geometry  not only generalizes 
the symplectic geometry of the cotangent bundle, 
but also provides information
about ``momentum frames'' along the par\-ti\-cle trajectories.
The equations also give rise to Poisson algebras of tensorial classical observables.

We find in the literature various linear${}^{\ref{Gu}-\ref{RR}}$ and
affine${}^{\ref{Ki2}-\ref{CCI}}$
models
of mul\-ti\-sym\-plec\-tic (mul\-ti\-var\-i\-able symplectic)
geometry on a covariant multiphase (multivariable phase) space of a
classical field. 
There are  two standard re\-pre\-sen\-ta\-tions  of 
the linear model of a multiphase space (Refs.\ \ref{Gu} and~\ref{KT}), but
Gotay${}^{\ref{Go}}$ has pointed out that both possess the
same inherent flaw.
Indeed, the linear model in either re\-pre\-sen\-ta\-tion has a covariant mul\-ti\-sym\-plec\-tic form,
 but such a form is not naturally
defined unless we arbitrarily  specify an Ehresmann connection on the
underlying field configuration space.
 The 
affine model suggested by Kijowski${}^{\ref{Ki2},\ref{Ki1}}$ and refined by 
Gotay, et al.,${}^{\ref{Go},\ref{GIMMsy}}$
is an improvement because 
the affine mul\-ti\-sym\-plec\-tic
structure can be
defined intrinsically.

In the affine model,  a variational principle applied to
a given field Lagrangian${}^{\ref{Go},\ref{GIMMsy}}$ 
would determine a unique Ehresmann connection. 
However, as the GIMMsy${}^{\ref{GIMMsy}}$  monograph demonstrates,
the problem with the affine model is that the space of
field momentum observables in the  affine multiphase space is not closed
under the naturally defined Poisson bracket.
The generalized symplectic geometry developed in this paper
addresses
both the problem with the connection in the linear model
and the problem with the bracket of momentum observables in the affine model.

In this paper we will show that
a particular symmetry breaking
of the linear frame bundle $LY$ of
the field configuration space $Y$ produces a subbundle, the 
{\em vertically adapted linear frame bundle} $\lvy$.
The $n$-symplectic geometry on $LY$ pulled back by inclusion to $\lvy$ 
generalizes both the  
linear and  affine mul\-ti\-sym\-plec\-tic geometries.  
We will  associate to $\lvy$ the multivelocity spaces {\em and} the
 multiphase spaces for both the linear and the affine models.
Then we can identify an Ehresmann connection on $Y$ with a symmetry breaking
of $\lvy$, 
thus addressing the shortcoming of the linear model.
To solve the problem in the affine model, we construct a momentum observable
from a projectable vector field on $Y$ and solve the generalized symplectic
structure equation on $\lvy$ to obtain a vector-valued Hamiltonian vector 
field on $\lvy$.
From the Hamiltonian vector fields we may define a Poisson bracket not only
that generalizes the Poisson
bracket on  the momentum observables in the affine model but also 
that makes the vector space of momentum observables into
an algebra under the bracket.

The format of this paper is as follows.
In Sec.\ II, we summarize both the linear and affine mul\-ti\-sym\-plec\-tic
geometries, and in Sec.\ III
we summarize $n$-symplectic geometry. 
We construct a
generalized symplectic geometry
on $\lvy$
 in Sec.\ IV\@. 
Secs.\ V  and VI
show that the multivelocity spaces and multiphase spaces
are associated to $\lvy$.
This work reveals  the one-dimensional vector bundle structure of the
affine 
multiphase space over the linear multiphase space.
Sec.\ VII is devoted to the generation of the affine 
mul\-ti\-sym\-plec\-tic potential
from the generalized symplectic potential of the vertically adapted
frame bundle. This leads to the resolution of the
problem with momentum observables.

\section{\normalsize\bf  MULTISYMPLECTIC GEOMETRY}
\label{sec:over}

\setcounter{equation}{0}

This section summarizes the linear and the affine models of mul\-ti\-sym\-plec\-tic geometry.
For an extensive bibliography see GIMMsy.${}^{\ref{GIMMsy}}$
Let $X$ be an oriented $n$-dimensional manifold and
let ${\pi_{XY}}: Y \rightarrow X$ be a fiber bundle
 with a $k$-dimensional fiber.  
(Note:  In general, $\pi_{BA}$ will denote a projection from $A$ onto $B$.)
A classical field is a section of 
the {\em field  configuration space}
$Y$ over the {\em parameter space} $X$.
From local coordinates $\{x^i\}, i = 1,\dots,n$, on $X$ we may construct local adapted coordinates $\{ x^i, y^A\}, i = 1,\dots,n, A = 1,\dots,k$ on $Y$.
Define the {\em vertical subbundle} of $TY$ to be 
$$V(TY) := \{ w_y\;|\; y\in Y , \, w_y \in T_yY \, \mbox{and}
\; \pi_{XY\,*}(w_y)=0 \}\, ,$$  
and define the {\em linear multivelocity space} to be $\hom_Y(TX,V(TY))$, the vector 
bundle over $Y$ whose fiber over $y\in Y$ is the vector space 
$\lin(T_{{\pi_{XY}}(y)}X,V(T_yY))$.
Define the {\em linear multiphase space} $J^*Y$
to be the vector bundle dual to the linear multivelocity space.
We shall discuss two equivalent re\-pre\-sen\-ta\-tions of $J^*Y$,
which we shall denote the G\"unther re\-pre\-sen\-ta\-tion${}^{\ref{Gu}}$
and the Kijowski and Tulczyjew (KT) re\-pre\-sen\-ta\-tion.${}^{\ref{KT}}$

Each fiber of the vector bundle 
$\hom_Y(TX,V(TY))^*$
is the  linear dual to the corresponding fiber of $\hom_Y(TX,V(TY))$.
Using the vector bundle isomorphism${}^{\ref{Gu}}$ 
$$
\hom_Y(V(TY),TX)  \rightarrow  \hom_Y(TX,V(TY))^*     :
\alpha_y    \mapsto \mbox{tr}(\alpha_y\circ \cdot \,) \, ,
$$
we obtain $\gunther$,
the G\"unther re\-pre\-sen\-ta\-tion of $J^*Y$. 
From  local adapted coordinates  on $Y$, 
define local canonical coordinates $\{x^i,y^A,p^j_B\}$
on $\gunther$
where 
$p^j_B(\alpha_y) := dx^j ( \alpha_y( \basisy B|_y ))$.
There does not exist a 
natural mul\-ti\-sym\-plec\-tic potential on $\gunther$ that is analogous to
the symplectic potential on the cotangent bundle.
Indeed,  $\alpha_y\in \gunther$ is defined only on 
vertical vectors in $T_yY$. 
Choosing an Ehresmann connection on $Y$ will remedy this 
problem.${}^{\ref{Go}}$ 
Let $\gamma$ be an Ehresmann connection on $Y$, represented as
 a projection $\gamma_y:T_yY\rightarrow V(T_yY)$ fibered over $Y$.
Now define a $TX$-valued 
one-form $\Theta_\gamma$ on $J^*Y = \gunther$ by 
$$
\Theta_\gamma (\alpha_y) = \pi^*_{Y\, J^*Y}(\alpha_y\circ\gamma)\, , \qquad \alpha_y\in J^*Y\; .
$$
G\"unther${}^{\ref{Gu}}$ calls $\Theta_\gamma$ a {\em canonical polysymplectic form},  although its definition requires an object extrinsic to the
natural geometry of $J^*Y$.
We shall call $\Theta_\gamma$  a {\em linear mul\-ti\-sym\-plec\-tic potential}
on $J^*Y$ in the G\"unther re\-pre\-sen\-ta\-tion   
because 
$\Theta_\gamma$ generates a {\em linear} mul\-ti\-sym\-plec\-tic geometry. 
In local coordinates, using the summation convention, 
$$
\Theta_\gamma  = \left(p^i_Ady^A  + p^i_A\gamma^A_j dx^j\right)\otimes
{\basisx i  } \; .
$$ 

The KT re\-pre\-sen\-ta\-tion of $J^*Y$ is 
$\hom_Y(V(TY),\wedge^{n-1}X)$ where $\wedge^{n-1}X$ is the bundle of 
$(n-1)$-forms on $X$,
and the values of sections of $\kt$  are 
vector-densities.${}^{\ref{KT}}\,$
In local coordinates  on  $X$ define
\beqa
d^nx & := & dx^1 \wedge \cdots \wedge dx^n \, , \nonumber \\   
d^{n-1} x_i  &: = & \basisx i \hook d^nx \, ,  \label{eqn:dn}\\
d^{n-2} x_{ij}  & : = & \basisx j \hook (\basisx i \hook d^nx)\, ,\quad \mbox{etc.,} \nonumber 
\eeqa
where $\hook$ denotes the inner product of a vector with a differential form.
Local coordinates on $\kt$ are $\{x^i, y^A, p^j_B\}$, where 
$$
p^j_B(\alpha_y) = 
(-1)^{j-1}\basisx n \hook \cdots \hook\widehat{\basisx j}
\hook \cdots \hook \basisx 1 \hook
\left( \alpha_y(\basisy B )\right)
$$
and $\widehat \basisx j$ denotes omission of $\basisx j$.
So   the mul\-ti\-sym\-plec\-tic potential on $J^*Y$ in the
KT re\-pre\-sen\-ta\-tion in local coordinates is
$$
\Theta_\gamma = \left(p^i_Ady^A  + p^i_A\gamma^A_j dx^j\right)\otimes
d^{n-1}x_i \; .
$$
Given 
a  volume form $\omega$ on $X$, we may
use the vector space isomorphism  
$ T_xX \rightarrow\wedge^{n-1}_x X:v \mapsto v \hook\omega$
to define a vector bundle isomorphism from the G\"unther re\-pre\-sen\-ta\-tion
of $J^*Y$ to the KT re\-pre\-sen\-ta\-tion.
In either re\-pre\-sen\-ta\-tion
 $d\Theta_\gamma$ is nondegenerate in the sense that for a vector field $X$, 
$X\hook d\Theta_\gamma = 0 $ if and only if $X = 0$.
We conclude that
 $(J^*Y, d\Theta_\gamma)$ is a {\em connection-dependent linear
 mul\-ti\-sym\-plec\-tic manifold}.

Ragionieri and Ricci${}^{\ref{RR}}$
use a hybrid re\-pre\-sen\-ta\-tion of $J^*Y$
using $TX \otimes_Y V^*Y $, which is isomorphic to $\gunther$ as
a vector bundle,
but they assume a volume form on $X$
and obtain a ``generalized Liouville $n$-form'' on  $TX \otimes_Y V^*Y $.
This form is also dependent upon the choice of a connection.

To avoid requiring a connection,
we must reconsider the  multivelocity space.
Ragionieri and Ricci
claim that the appropriate multivelocity 
space is not $\hom_Y(TX,V(TY))$, but rather, the (first-order) {\em jet bundle} $JY$,
the  affine bundle over
$Y$ whose fiber over $y\in Y$ consists of linear maps${}^{\ref{GIMMsy}}$
$\gamma_y: T_{{\pi_{XY}}(y)}X\rightarrow T_yY$ satisfying
${\pi_{XY}}_*\circ \gamma_y = \mbox{Id}\,_{T_{{\pi_{XY}}(y)}X}$.
The affine structure of $JY$ comes from a difference function
$\delta : JY\times JY \rightarrow \hom_Y(TX,V(TY))$. 
Recall${}^{\ref{GIMMsy}}$ that a section of $JY$
over $Y$ can be identified
 with an Ehresmann connection on $\pi_{XY}:Y\rightarrow X$.

 The {\em  bundle of affine cojets}$\,{}^{\ref{Go},\ref{GIMMsy}}$ is the 
vector bundle  $J^\star Y$ over $Y$ whose fiber over $y\in Y$ is the set
of affine maps from $J_yY$ to $\wedge^{n}_{{\pi_{XY}}(y)}X.$
It follows that  
$\dim J^\star Y  = \dim JY + 1$.
We find an equivalent
description of $J^\star Y$ to be useful.
The {\em affine multiphase space} $Z$
is the bundle of $n$-forms on $Y$ whose 
fiber $Z_y$ 
over $y$ is 
$$Z_y := \{ z\in \wedge^n_yY  \, | \, v\hook w\hook z = 0 \; \forall v,\, w \in
 V(T_yY)\}\; .$$
The affine multiphase space, originally defined by Kijowski,${}^{\ref{Ki2}}$ 
admits an
$n$-form, 
$\Theta_z = \pi^*_{YZ}(z)$,
which is the pullback via inclusion of the canonical $n$-form on $\wedge^nY$.
(Kijowski${}^{\ref{Ki1}}$ attributes the notion of a multiphase space
 to an unpublished result of Tulczyjew.)
 GIMMsy${}^{\ref{GIMMsy}}$ proves that $Z$ is 
``canonically'' isomorphic to $J^\star Y$.
On $Z$, we can define local coordinates $\{ x^i, y^A, p^j_B, p\}$ 
where
\begin{eqnarray}
  p(z) & = &  \basisx n \hook\cdots\hook\basisx 1
\hook z \quad\mbox{and} \label{eqn:zcoords}\\
p^j_B(z)
& = & (-1)^{j-1} \, \basisx n \hook \basisx {n-1} \hook\dots\hook 
\widehat{\basisx j} \hook \dots\hook\basisx 1\hook\basisy B \hook z \, .
\nonumber \end{eqnarray}
So in local coordinates 
\beq\label{eqn:Theta}
\Theta = p^i_A dy^A \wedge d^{n-1}   x_{i} +  p d^nx
\eeq
where $d^nx$ and $d^{n-1}x_i$ given in~(\ref{eqn:dn}) are
pulled up from $X$ to $Z$.
The $(n+1)$-form $d\Theta$ is nondegenerate, so the pair $(Z,d\Theta)$ is an 
affine {\em  mul\-ti\-sym\-plec\-tic manifold}.${}^{\ref{GIMMsy}}$  
A connection is {\bf not} required to 
define the affine mul\-ti\-sym\-plec\-tic structure.

Let ${\cal X}Y$ be the Lie algebra of  vector fields on $Y$.
Denote the space of
vector 
fields of $Y$ projectable to $X$ by $\projectable$.
Note that $\projectable$ is a Lie subalgebra
of ${\cal X}Y$, since $[{\pi_{XY}}_*v,{\pi_{XY}}_*w] = {\pi_{XY}}_*[v,w]$.
(That is, the Lie brackets are $\pi_{XY}$-related.${}^{\ref{Ch}}$)

\define 
Let $v\in\projectable$. A {\em momentum observable}$\,{}^{\ref{GIMMsy}}$ based
on $v$ is an ${(n-1)}$-form $f_v$ on $Z$ defined by
$$f_v(z) := \pi_{YZ}^*(v\hook z)\, .$$
Let $T^1(Z)$ denote the vector space of momentum observables.

\bigskip

\noindent If in local adapted coordinates on $Y$, $v = v^i(x^j)\basisx i + v^A(x^j,y^B)\basisy A$, then in local coordinates on $Z$,
\beq\label{eqn:fv}
f_v (z)
 =  (p^i_Av^A + pv^i)d^{n-1}x_i - p^i_A v^j dy^A\wedge d^{n-2}x_{ij}\; .
\eeq
The {\em Hamiltonian vector field} $X_{f_v}$ is obtained from $f_v$ 
via the {\em mul\-ti\-sym\-plec\-tic structure equation}
\beq\label{eqn:multistruc}
df_v = - X_{f_v}\hook d\Theta \, .
\eeq
The local coordinate expression for $X_{f_v}$ is
\beq \label{eqn:xfvcoords}
X_{f_v}  = v^k \basisx k + v^A \basisy A 
+  \left( p^j_A \partialx {v^i} j - p^i_A\partialx {v^j} j
- p^i_B \partialy {v^B} A \right)   \basisp i A
-  \left( p\partialx {v^i} i  + p^i_A\partialx {v^A} i  \right)   
\basisp {} {}\, .
\eeq

From equation (\ref{eqn:xfvcoords}) and Ref.~\ref{Ch}, if 
$v,w \in \projectable$ then 
$\pi_{YZ*}\left[X_{f_v},X_{f_w} \right] =  [v,w]$.
From (\ref{eqn:Theta}), (\ref{eqn:fv}) and (\ref{eqn:xfvcoords}),
\beq\label{eqn:xhooktheta}
X_{f_v}\hook \Theta  =    f_v\; .
\eeq
Using (\ref{eqn:xhooktheta}),
\beq\label{eqn:brackhooktheta}
\left[X_{f_v},X_{f_w} \right] \hook\Theta_z 
 =  \pi^*_{YZ}\left( [v,w]\hook z\right) 
 =  f_{[v,w]}(z)\; . 
\eeq

\noindent
Using equations~(\ref{eqn:multistruc}) and~(\ref{eqn:xhooktheta})
and the Lie derivative identity${}^{\ref{AM}}$ 
\beq\label{eqn:liederiv}
{\cal L}_X\alpha = X \hook d\alpha + d(X\hook \alpha)
\eeq
it follows that
\beq\label{eqn:lieconstant}
{\cal L}_{X_{f_v}}\Theta   =   0 \; .
\eeq

\define Let $f_v$ and $ f_w \in T^1(Z)$ and let $X_{f_v}$ and $X_{f_w}$ be their corresponding
Hamiltonian vector fields, respectively.  
Define the {\em Poisson bracket}$\,{}^{\ref{GIMMsy}}$
 of $f_v$ and $f_w$  to be
$$
\{f_v,f_w\} : =   - X_{f_v}\hook (X_{f_w} \hook d\Theta).
$$

\noindent The Poisson bracket is not a true Poisson bracket 
because there lacks an 
associative multiplication of $(n-1)$--forms on which the bracket acts
as a derivation.
Using ${\cal L}_{X_{f_v}}\Theta   =   0$,
 equations (\ref{eqn:brackhooktheta}), (\ref{eqn:liederiv}), (\ref{eqn:lieconstant}), 
and the identity${}^{\ref{AM}}$
$$
[X,Y]\hook \alpha = {\cal L}_X(Y\hook \alpha) - Y\hook({\cal L}_X\alpha)\; ,
$$
we obtain
\beq\label{eqn:pbexact}
\{f_v,f_w\} 
=  f_{[v,w]}
        - d(X_{f_v}\hook(X_{f_w}\hook\Theta)) \; .  
\eeq
The exact form on the right side of equation (\ref{eqn:pbexact}) is not in $T^1(Z)$
because from equation~(\ref{eqn:multistruc})
the Hamiltonian vector field of an exact form is the zero vector field on $Z$, but
the momentum observable corresponding to the zero vector field is
the zero $(n-1)$-form on $Z$.
Thus, the space  $T^1(Z)$
is {\bf not} closed under the Poisson bracket.
Equation (\ref{eqn:pbexact}) explains the remark in GIMMsy${}^{\ref{GIMMsy}}$
 that the
Poisson bracket of two momentum observables ``is up to the addition of
exact terms, another momentum observable.''

\bigskip

\noindent{\bf Definitions} \begin{description}

\item[{$\bullet$}] The Lie algebra of  {\em locally Hamiltonian vector fields} on $Z$
 is 
$$LHV^1(Z) := \{X\in\,{\cal X} Z \;|\; {\cal L}_X d\Theta = 0 \}.$$

\item[{$\bullet$}] The vector space of {\em allowable Hamiltonian observables} on $Z$
is 
$$
HF^1(Z) :=
 \{ f\in\, \wedge^{n-1}Z \;|\; d f = -X\hook d\Theta\, , \,X \in LHV^1(Z)\}\, .
$$

\item[{$\bullet$}]The vector space of  {\em Hamiltonian vector fields} on $Z$
 is   
$$HV^1(Z) := \{X\in\,{\cal X} Z  \;|\;X\hook d\Theta = -df , \, 
f \in HF^1(Z) \}\, .$$

\end{description}

By equation (\ref{eqn:multistruc}), a momentum observable on $Z$ is an allowable Hamiltonian observable.
Define the {\em Euler vector field}$\,{}^{\ref{Wo}}$ ${\cal E}$ on $Z$ by
${\cal E}\hook d\Theta = \Theta$.
Because $d\Theta$ is nondegenerate, ${\cal E}$ is well-defined.
If $f\in T^1(Z)$, then, 
using (\ref{eqn:multistruc}) and (\ref{eqn:xhooktheta}),
\beq\label{eqn:Euler}
{\cal E}\hook df  = {\cal E}\hook(-X_{f}\hook d \Theta)
 =   X_f\hook \Theta 
 =  f \; .
\eeq
Define a linear operator $\tilde {\cal E}$ on $\wedge^{n-1}Z$ by 
$\tilde {\cal E} (f) := {\cal E} \hook df\; .$
When restricted to  $HF^1(Z)$,
$\tilde {\cal E}$ is a projection operator.
If $f\in T^1(Z)$ then $\tilde{\cal E}(f )= f$
by~(\ref{eqn:Euler}).
If $n=2$ and $k=2$ or
if $n\geq 3$ then
$HF^1(Z)$ is an algebra
under the Poisson bracket extended to $HF^1(Z)$ and  
the image of $HF^1(Z)$ under $\tilde{\cal E}$ is $T^1(Z)$.  Thus,
$$HF^1(Z) = T^1(Z) \oplus \mbox{Ker}\, (\tilde{\cal E}| HF^1(Z)).$$ 
Furthermore, the nonzero exact $(n-1)$-forms on $Z$ are in 
$\mbox{Ker}\, (\tilde{\cal E}| HF^1(Z))$, so  $T^1(Z)$ is a vector 
subspace of $HF^1(Z)$, but  by
equation~(\ref{eqn:pbexact}), $T^1(Z)$ is not a subalgebra of $HF^1(Z)$.
Define a re\-pre\-sen\-ta\-tion $HF^1(Z)\rightarrow HV^1(Z): f\mapsto X_f$ and 
it follows that
\beq\label{eqn:HVtoHFZ}
HV^1(Z) \simeq HF^1(Z)/ \{\mbox{closed $(n-1)$--forms}\} \; .
\eeq

If $n=1$ then $Z=T^*Y$, $HF^1(Z) = C^\infty(T^*Y)$, the Poisson
bracket is the usual Poisson bracket on the cotangent bundle $T^*Y$, and $HV^1(Z) = HF^1(Z)/\r$.

\section{\normalsize\bf  n-SYMPLECTIC GEOMETRY ON THE LINEAR FRAME BUNDLE}
\label{sec:norris}
\setcounter{equation}{0}

This section is a summary
of Norris's theory  of $n$-symplectic geometry on 
the linear frame bundle.${}^{\ref{No1},\ref{No2}}$
Let $M$ be an $n$-dimensional manifold and  let $\tau:LM\rightarrow M$ be  the
 bundle of linear frames of $M$. That is, 
$$LM := \{(x, e_i) \, | \, x\in M , \, \{e_i\}  \,
\mbox{is a frame of}\; T_xM\}$$  
The structure group of $LM$ is the real general linear group 
$GL(n)$, which acts freely
 on the right of $LM$.
Also, $LM$  supports an
 ${\r}^n$-valued  {\em soldering  one-form${}^{\ref{KN}}$ 
$\theta$},
 defined by
$$
\theta(Y){:=} u^{-1}(\tau_* Y)\ \ \ \ \forall\  Y\in T_uLM\ ,
$$
where if  $\{r_i\}, i=1,2,\dots,n$, denotes the
standard  basis of $\r^n$, then
$u=(x,e_i):{\r}^n    \rightarrow  T_{\tau(u)}M$ is the linear isomorphism
 $\xi^i r_i  \mapsto  \xi^ie_i$. 
Compare $\theta = \theta^i r_i$
to the real-valued  canonical one-form $\vartheta$ on
$T^*M$.
The ${\r}^n$-valued two-form $d\theta$ 
 on $LM$ 
is closed and nondegenerate,
so $d\theta$ 
is considered to be
an $\r^n$-valued symplectic structure,
or an {\em $n$-symplectic structure} on $LM$.

The theory of $n$-symplectic geometry on $(LM,d\theta)$  is
based on the
generalized structure equation
\beq\label{eqn2:struc1}
d\hat f^{i_1i_2\dots i_p}=-p!X_{\hat f}^{i_1i_2\dots i_{p-1}}\hook d\theta^{i_p}\, ,
\eeq 
where the functions  $\hat f^{i_1i_2\dots i_p}$ are the components of a $\otimes^p\r^n$-valued function  $\hat f$,
and the vector fields $X_{\hat f}^{i_1i_2\dots i_{p-1}}$ are the components of a
$\otimes^{p-1}\r^n$-valued vector field $ X_{\hat f}$.
Because  for $g\in GL(n)$ the {\em $n$-symplectic potential}
$\theta$  transforms tensorially under right translations $R_g$
 by $R^*_g\theta=g^{-1}\cdot\theta$,
not every
 $\otimes^p\r^n$-valued function  $\hat f$ is  compatible with 
equation (\ref{eqn2:struc1}). Hence $n$-symplectic geometry selects classes
of allowable observables,  in contrast with the fact that
all smooth real-valued functions on $T^*M$
are allowable observables.

The allowable observables in $n$-symplectic geometry divide naturally into the
symmetric and antisymmetric Hamiltonian functions,  which we denote
by $SHF$ and $AHF$, respectively. The space $SHF$ is the direct sum
$\oplus_{p=1}^\infty SHF^p$ where $SHF^p$ is the space of $(\otimes_s)^p\r^n$-valued
functions defined on $LM$ 
 that are compatible with a symmetrized version of 
(\ref{eqn2:struc1}) and  $\otimes_s$ denotes
the symmetric tensor product.  The elements of $SHF^p$ are
 degree $p$ polynomials in the generalized local momentum coordinates on
$LM$ with coefficients that are constant on the fibers of $LM$.
 In particular,
the elements of $SHF^p$ whose local coordinate representatives are homogeneous
polynomials
are in bijective correspondence with symmetric degree $p$
contravariant tensor fields on $M$ and are precisely
the tensorial elements of $SHF^p$.
We call this subset $ST^p$ (for {\bf S}ymmetric and {\bf T}ensorial) and 
define $ST=\oplus_{p=1}^\infty ST^p{}^{\ref{No1}}$.

For $ST^p$, $p\geq 1$, equation  (\ref{eqn2:struc1}) is replaced by
\beq\label{eqn:struc2}
d\hat f^{i_1i_2\dots i_p}=-p!X_{\hat f}^{(i_1i_2\dots i_{p-1}}\hook d\theta^{i_p)}
\eeq
where the parentheses  denote symmetrization over the enclosed indices.
For each $\hat f\in ST^p$ this equation determines an equivalence class
 of $({}^{n+p-2}_{\ \  p-1})$
 vector fields.
 If $p=1$ then the equivalence relation is equality. If  
$p\geq 2$ then the symmetrization on indices in
(\ref{eqn:struc2}) introduces a degeneracy which is not present for $p=1$. 
More precisely, for each $\hat f\in ST^p$ and each $u\in LM$ there exists a
local
section $\xhatf=\xhatf^{\multii {p-1}}(r_{i_1}\otimes_s\cdots\otimes_s  r_{i_{p-1}})$
of the vector bundle $T(LM)\otimes(\otimes_s)^{p-1}\r^n\rightarrow LM$
defined on an open subset 
$U_u$ of $LM$ such that the  vector fields
$\xhatf^{\multii {p-1}}$ satisfy (\ref{eqn:struc2}) for all $(i_1,i_2,\dots,i_{p})$.
In fact there exists more than one such local section since
$$
K_v=\{Y_v\ |\ Y_v^{(i_1\dots i_{p-1}}\hook d\theta^{i_p)}=0\ \hbox{for all}\ (i_1,i_2,\dots,i_{p})\}
$$
is a nontrivial subspace of $T_v(LM)\otimes(\otimes_s)^{p-1}\r^n$ for each
$v\in LM$.  Indeed, 
$K=\bigcup_{v\in LM}(K_v)$
 is a vector subbundle of $T(LM)\otimes(\otimes_s)^{p-1}\r^n$
and for each $\hat f\in ST^p$ there exists a unique global section $\sigma$ of
$(T(LM)\otimes(\otimes_s)^{p-1}\r^n)/K\rightarrow LM$ such that for
each $u\in LM$ there exists a
neighborhood $U$ of $u$ and a local section $\xhatf$ of 
$T(U)\otimes(\otimes_s)^{p-1}\r^n\rightarrow U$
having the property that $\xhatf$ satisfies (\ref{eqn:struc2}) and
$\sigma(v)=\xhatf(v)+K_v$ for each $v\in U$.  Denote $\sigma$ by
$\vectorclass f$.  If we fix $I=(i_1,\dots, i_{p-1})$ then there is also a
subbundle $K^I$ of $T(LM)$ such that for $v\in LM$,
$$
K^I_v=\{\ Y^{\multii {p-1}}_v\ |\ Y^{(\multii{p-1}}_v\hook d\theta^{i_p)}=0\ 
\hbox{for all}\ i_p\ \}\, .
$$
Moreover there is a unique section $\sigma^I$ of $T(LM)/K^I\rightarrow LM$
such that for  $u\in LM$, 
there is a neighborhood $U$ of $u$ and a local section 
$\xhatf^i$ of $T(U) \rightarrow U$ with the property that
$\xhatf^I$ satisfies equation (\ref{eqn:struc2})
 and
$\sigma^I(v)=\xhatf^I(v)+K_v^I$ for all
$v \in U$.
We denote this section $\sigma^I$ by $\vectorclass f^I
=\vectorclass f^{\multii {p-1}}$.

The fact that elements of $ST^p$  determine
 equivalence classes of vector fields
 does not affect
the basic algebraic structures in $n$-symplectic geometry. 
For each $p\geq 1$ the set of
equivalence classes of $\otimes^{p-1}_s{\r}^n$-valued vector fields on $LM$  forms an
 infinite-dimensional vector space.  Denote by $HV(ST^p)$ the
vector space of $\otimes_s^{p-1}{\r}^n$-valued equivalence classes of vector fields determined 
 by elements of $ST^p$ by equation~(\ref{eqn:struc2}).
For $\hat f\in ST^p$ and $\hat g\in ST^q$ define the Poisson bracket 
$\{\ ,\ \}:ST^p\times ST^q\to ST^{p+q-1}$
by
\beq\label{eqn:Poissonn}
\pb f g^{i_1 i_2\dots i_{p+q-1}} =p!\xhat f^{(i_1 i_2 \dots i_{p-1}}\left(\hat g^{i_p i_{p+1} \dots
i_{p+q-1})}\right)
\eeq
where $\xhat f^{i_1 i_2 \dots i_{p-1}}$ is any representative of the
equivalence class $\vectorclass f^{i_1 i_2 \dots i_{p-1}}$.  The bracket 
defined in (\ref{eqn:Poissonn}) is easily shown to be independent
of the choice of representatives and it has all the properties of a Poisson
bracket.  In fact  when  the bracket in (\ref{eqn:Poissonn})  is
re-expressed on the base manifold $M$, it gives the differential concomitant of Schouten${}^{\ref{Sc}}$ and Nijenhuis${}^{\ref{Ni}}$ of the symmetric tensor fields corresponding to $\hat f$ and $\hat g$.

\begin{thm} (Norris$\,{}^{\ref{No1}}$) The space $ST$ of symmetric tensorial
functions on $LM$ is a Poisson algebra with respect to the Poisson bracket defined in~(\ref{eqn:Poissonn}).
\end{thm}

Denote the direct sum of the vector spaces $HV(ST^p)$ by $HV(ST)$. 
There is a well-defined Lie bracket
on the vector space $HV(ST)$ that satisfies
$$
[\vectorclass f,\vectorclass g]=\lbrack\!\lbrack
X_{\pb f g}\rbrack\!\rbrack\, .
$$

A degree one
 tensorial observable
$\hat f\in T^1(LM)= ST^1$ corresponds to a unique vector field
${\boldmath f}$ on $M$.${}^{\ref{KN}}$ 
 In local canonical coordinates
$\{x^i, \pi^j_k := e^j(\basisx k)\}$ 
on $LM$ we may write
$\hat f^i= f^j(x)\pi_j^i$ where
${\boldmath f}=f^i(x){\partial\over\partial x^i}$, and
$x\in M$. 
For $p=1$ equation (\ref{eqn:struc2}) now 
has a unique solution $X_{\hat f}\in HV(ST^1)$,
given in local coordinates by
$$
X_{\hat f}=f^i(x){\partial\over\partial x^i}-{\partial f^i\over\partial x^j}\pi^k_i
{\partial\over\partial \pi^k_j}\, .
$$
The vector field $X_{\hat f}$ is 
the {\it natural lift}$\,{}^{\ref{KN}}$
of the vector field ${\boldmath f}$ to $LM$.
See Ref.~\ref{No1} for  discussions of the full Poisson algebra $SHF$
and of the graded Poisson algebra $AHF$.

All of the  basic features of symplectic geometry on $T^*M$ are induced from the $n$-symplectic geometry on $LM$.
Indeed, it is well known that $T^*M$ is isomorphic to the associated bundle
 $LM\times_{GL(n)} \r^{n*}$. Furthermore, 
the relationship
between the canonical one-form $\vartheta$ on $T^*M$ and the soldering one-form $\theta$ on $LM$ is${}^{\ref{No2}}$
\beq\label{eqn:Sniatycki}
\vartheta_{[(u,\alpha)]}(\tilde X)= \,<\theta_u(X),\alpha>
\eeq
where $u \in LM$, $[u,\alpha] \in T^*M \simeq LM\times_{GL(n)}\r^{n*}$
and
$\tilde X$ is a tangent
vector at $[u,\alpha]$ that projects to the same vector as the tangent 
vector $X$ at $u$.  
Thus the symplectic potential $\vartheta$ for
symplectic geometry on $T^*M$ is induced from the soldering one-form
$\theta$ on $LM$. 
Moreover, if for some nonzero $\alpha\in \r^{n*}$ we define a map${}^{\ref{No2}}$
$$
\phi_\alpha : LM\rightarrow T^*M : u\mapsto [u,\alpha] \, ,
$$
then 
the range of $\phi_\alpha$ is $T^*M$ with the zero-section
deleted, and thus (\ref{eqn:Sniatycki}) becomes
\beq\label{eqn:one-form}
\phi^*_\alpha \vartheta = \, <\theta,\alpha>\, .
\eeq
Also 
the homogeneous degree $p$ polynomial observables on \mbox{$T^*M$} are induced from 
 elements of $ST^p$. 
Indeed,
 for $\hat f\in ST^p$ and $u \in  LM$, define 
\beq\label{eqn:alphas}
\tilde f:T^*M\to\r \, : \, [u,\alpha] \mapsto \,<\hat f(u),\overbrace{(\alpha,\dots,\alpha)}^p >\, .
\eeq
Because 
$\hat f$ is tensorial, $\tilde f$ is well defined.
In local coordinates, $\pi^i_j(x,e_k)\alpha_i=e^i(\basisx j)\alpha_i=
p_j(e^i\alpha_i)$ where $\{p_j\}$ are the standard local momentum coordinates
on \mbox{$T^*M$}.
For example, if $p=2$ and  
$\hat f= \hat f^{kl}(x) \pi^i_k \pi^j_l r_i\otimes_s r_j \in ST^2$ then 
$$
\tilde f([u,\alpha])=f^{kl}(x)p_kp_l \,  .
$$

\begin{thm}\label{thm:lmvecfields} (Norris$\,{}^{\ref{No2}}$)
Let $\hat f\in ST^p$,  let $\vectorclass f$ be
the associated equivalence class of Hamiltonian vector fields determined 
by~(\ref{eqn:struc2}), 
and let $\tilde f$ be the degree p homogeneous polynomial observable on \mbox{$T^*M$} determined by $\hat f$ as in~(\ref{eqn:alphas}). Then 
$$
X=p!\phi_{\alpha*}(\xhat f^{i_1i_2\dots i_{p-1}}\alpha_{i_1}\alpha_{i_2}\cdots
\alpha_{i_{p-1}})
$$
where  $\xhat f^{i_1i_2\dots i_{p-1}}$ denotes any set of representatives
of $\vectorclass f$,
$X$ is a vector field on $T^*M$
with the zero-section of $T^*M$ deleted,
 and $X = X_{\tilde f}$.
\end{thm}

\section{\normalsize\bf  THE VERTICALLY ADAPTED LINEAR FRAME BUNDLE}\label{sec:valfb}
\setcounter{equation}{0}

We now present
a new principal bundle gen\-e\-ral\-i\-za\-tion
of mul\-ti\-sym\-plec\-tic geometry. 
From Norris's theory applied to the 
$(n+k)$-dimensional field  configuration space
$\pi_{XY}:Y\rightarrow X$,
the linear frame bundle $LY$ over $Y$ 
has an $(n+k)$-symplectic geometry. 

\define The {\em vertically adapted frame bundle} $\lvy$ 
is defined as
$$
\lvy  := \{\ptlvy \in LY \;|\; 
\{\epsilon_A\} \hbox{ is a frame of } V(T_yY) \}\, .
$$
The terminology is motivated by the definition of an {\em adapted frame}.${}^{\ref{KN}}\,$  Again, as in 
Sec.~II, 
 $i = 1,\dots,n$ and $A = 1,\dots,k$.

\medskip

Let $M_k(\r^{n+k})$ denote the
$k^{th}$ Grassmann manifold of  $\r^{n+k}$.
A point in $M_k(\rnplusk)$ 
is the span of a set of $k$ linearly independent
 vectors $\{w_A\}$ in $\rnplusk$.
Define a left
action of $GL(n+k)$ on $M_k(\rnplusk)$
by 
\beq\label{eqn:grass}
 g\cdot \mbox{span}\{w_A\} := \mbox{span}\{g^B_A w_B\} \, .
\eeq

Let $i_1: \r^n\rightarrow\r^{n+k}: v\mapsto \hat v $
be the inclusion into the first $n$ slots of $\rnplusk$
and let $i_2: \r^k\rightarrow\r^{n+k}: w\mapsto \hat w $
be the inclusion into the last $k$ slots.
Let $\{r_i\}$ be the standard basis of $\r^n$ and
let $\{s_A\}$ be the standard basis of $\r^k$.
Then the isotropy
subgroup of $\mbox{span}\{\hat s_A \}$ expressed in a matrix re\-pre\-sen\-ta\-tion with respect to
the standard basis of $\rnplusk$ is the {\em adapted linear group},
$$
\ga  := \left\{ \left( \begin{array}{cc} 
N&0\\ A&K 
\end{array} 
\right)   \left. \;|\; N \in GL(n),\, K \in GL(k),\, A \in 
{\rkn} \right. \right\} \, .
$$
The group $\ga $ is a semidirect product of $GL(n)\times GL(k)$
and $\rkn$. For convenience we write $\nka\in \ga$ as $(N,K,A)$.
Since the action of $GL(n+k)$ on $M_k(\rnplusk)$ is transitive, 
we have proven that 
$M_k(\rnplusk)$ is diffeomorphic to $GL(n+k) / \ga $.

\begin{thm}
Define the map 
$$\phi:  LY \rightarrow M_k(\rnplusk): 
 w \mapsto w^{-1}(V(T_{\pi_{Y\,LY}(w)}Y))
$$
where $w\in LY $ is a linear isomorphism 
$w: \rnplusk \rightarrow T_{\pi_{Y\, LY}(w)}Y^{\ref{KN}}$.
  Then
$\phi$ is a symmetry-breaking map and $\lvy = \phi^{-1}(\mbox{\em span}\{\hat s_A\})$.
Thus, $\lvy$ is a symmetry-broken subbundle of $LY$ with structure group $\ga$.
\end{thm}

\noindent{\bf Proof}~
The action in (\ref{eqn:grass}) is transitive 
and $\phi$ is tensorial. 
Indeed,
$$\phi(w\cdot g) = (w\cdot g)^{-1}(V(T_{\pi_{Y\,LY}(w)}Y)) = g^{-1}(w^{-1}(V(T_{\pi_{Y\,LY}(w)}Y))) = g^{-1}(\phi(w)).$$
So, $\phi$ is a symmetry-breaking map.${}^{\ref{Tr}}$
Finally, if $w = (y,E_\mu) \in \lvy$, $\mu = 1,\dots, n+k$, then
$\phi(w)= \mbox{span}\{\hat s_A\}$  if and only if $V(T_yY) = w(\mbox{span}\{\hat s_A\} )
 := \mbox{span}\{  w(\hat s_A)\}$.
But $w(\hat s_A) = (y, E_\mu)(\hat s_A)   = E_{n+A}$.
So $w= (y, E_\mu)\in \phi^{-1}(\mbox{span}\{\hat s_A\}) $
if and only if $V(T_yY) = \mbox{span}\{ E_{n+A} \}_{A=1}^k$
if and only if $w  \in \lvy$.
Thus, the structure group of $\lvy$ is the isotropy subgroup of
$\mbox{span}\{\hat s_A\}$.~\blob

\bigskip

\define The {\em bundle of vertical frames} of $Y$, denoted $\vf$, is defined by
$$\vf := \left\{ (y,\epsilon_A)\, | \, y\in Y, \{\epsilon_A\} \hbox{ is a frame of } V(T_yY) \right\}.$$ 

\bigskip

\noindent The bundle $\vf\rightarrow Y$ is a principal fiber bundle  with structure group
 $GL(k)$.
 Let 
$\tau : LX \rightarrow X$ be the linear frame bundle over $X$.
Using $\tau$ we construct first the 
pullback bundle   
$\tau^* Y$ over $X$ and then
the pullback bundle $ \tauvf$ over $Y$.  See diagram~(\ref{diag:pullback}).
The elements of $\tauvf$ are of the form $\pttauvf$ where 
$(\pi_{XY}(y),f_i) \in LX$
and $(y, \epsilon_A) \in \vf$.
\beq\label{diag:pullback}
\begin{array}{ccc}
&&\\
\tauvf&\longrightarrow&\vf\\
&&\\
\downarrow&&\downarrow\\
&&\\
\tau^*{Y}&{\longrightarrow}&Y\\
&&\\
\downarrow&&\downarrow{\pi_{XY}}\\
&&\\
LX&{\tau\atop \longrightarrow}&X
\end{array}
\eeq

Represent $GL(n)\times GL(k)$  as a subgroup of 
$GL(n+k)$, identifying the ordered pair $(N,K)$ with the matrix
$ \left( \begin{array}{cc} N&0\\ 0&K \end{array} \right)$.
Observe that  $GL(n) \times GL(k)$  is a subgroup of $\ga$.
Define the right action of $GL(n) \times GL(k)$ on $\tauvf$ by
\beq\label{eqn:tauvfaction}\pttauvf\cdot(N,K) = \{(\pi_{XY}(y),y, f_jN^j_i, \epsilon_BK^B_A)\}\; .
\eeq
The action in~(\ref{eqn:tauvfaction}) is free and we see from
 diagram~(\ref{diag:pullback})
that $\tauvf$ is locally trivial over $Y$.
Thus $\pi_\tau: \tauvf \rightarrow Y$
 is a principal fiber bundle with structure group $GL(n)\times GL(k)$.
 The bundle $\lvy$ projects smoothly to
$\tauvf$ by the map 
$$\eta: \lvy \rightarrow \tauvf : (y,\frame)\mapsto ({\pi_{XY}}(y),y,{\pi_{XY}}_*(e_i),\epsilon_A).$$

\begin{lemma}\label{lemma:pfb2}
The bundle $\eta: \lvy \rightarrow \tauvf$ is a principal fiber bundle with structure group $(\r^{k\times n}, + )$.
Furthermore, the projection map $\eta$ is $GL(n)\times GL(k)$-equivariant.
\end{lemma}

\noindent{\bf Proof } 
The right action of $\ga $ on $\lvy$ is given by 
\beq\label{eqn:galvy}
\ptlvy \cdot (N,K,A) := (y,\{e_jN^j_i+\epsilon_BA^B_i,\epsilon_BK^B_A\})
\eeq
and this action is clearly free.
Any two points in $\eta^{-1}(x,y,f_i,\epsilon_A)$ are 
related by a unique translation by an element of 
 $(\r^{k\times n},+)$ expressed as the subgroup 
$\{(I_n,I_k,W), W \in \rkn\}$
of $\ga$. 
 Indeed, if 
$\eta(y,\{e_i,\epsilon_A\})=\eta(y',\{e'_i,\epsilon'_A\})$ then $y=y', \; 
\epsilon_A = \epsilon'_A$, and ${\pi_{XY}}_*e_i' = {\pi_{XY}}_*e_i$.
Thus
 the unique 
element of $\ga$ that maps $\ptlvy$ to $(y,\{e_i',\epsilon_A'\})$ is
$(I_n,I_k,A)$, so $e_i' = e_i + A^A_i\epsilon_A$.
Also, a local trivialization 
$\eta^{-1}(\pi_\tau^{-1}(U)) \simeq \pi_\tau^{-1}(U)\times
 \rkn$ can be obtained
in some open neighborhood $U$ of $y\in Y$.
Finally, using (\ref{eqn:tauvfaction}) and (\ref{eqn:galvy}),
\beq\label{eqn:etaequiv}
\eta(\ptlvy\cdot (N,K,0))  = \eta \ptlvy\cdot(N,K) \;.~\blob
\eeq

\bigskip

Let $(y,\{e^i, \epsilon^A\})$ be the coframe dual to the vertically
adapted frame $\ptlvy$.
We may define the right action of $\ga $ on 
$(y,\{e^i,\epsilon^A\})$ by   
\beq\label{eqn:coframeaction}
(y,\{e^i,\epsilon^A\})\cdot(N,K,A)
 =  (y,\{(N^{-1})^i_je^j,-(K^{-1}AN^{-1})^A_je^j  + (K^{-1})^A_B\epsilon^B\})\, . \eeq
The  coframe in~(\ref{eqn:coframeaction}) is dual to the frame
$\ptlvy\cdot(N,K,A)$.

Define a left action of $\ga$ on $\rkn$ by
\beq\label{eqn:gaaction}
(N,K,A)\cdot W := KWN^{-1}- AN^{-1}\;.
\eeq
The action 
in (\ref{eqn:gaaction}) transitive 
and the isotropy subgroup of $\{0\}$ is $GL(n)\times GL(k)$.

\begin{lemma}\label{lemma:taufvrkn}
$\tauvf\times\rkn\rightarrow Y$  is a principal fiber bundle with
structure group $\ga $.  
\end{lemma}

\proof First, define a right action of $\ga$ on $\tauvf\times\rkn$ by
\beq\label{eqn:tauvfrkn}
(u,W)\cdot(N,K,A) = \left(u\cdot(N,K), (N,K,A)^{-1}\cdot W\right) \, ,
\eeq
 where the action
on $\tauvf$ is given by~(\ref{eqn:tauvfaction}) and the $\ga $-action on  $\rkn$
is given by (\ref{eqn:gaaction}).
If $(u,W)\cdot(N,K,A) = (u,W)$  then  $u\cdot(N,K) = u$
and $(N,K,A)^{-1}\cdot W = W$.  The
action on $\tauvf$ is free, so $N=I_n$ and $K=I_k$,
and thus  $(I_n, I_k,A)^{-1}\cdot W = W$,
so $A=0$.  Thus the $\ga $-action is free.
The projection from
$\tauvf\times\rkn$ onto $Y$  is $\pi_\tau\circ \Pr_1$,
where $\Pr_1$ is projection in the first slot of the Cartesian product.
So for $y\in Y$ there exists an open neighborhood $U$ of $y$ such that
$\Pr_1^{-1}(\pi_\tau^{-1}(U)) \simeq  U \times \ga$.~\blob

\begin{thm}\label{thm:five}  The following are in pairwise bijective correspondence: 
\begin{list}
{(\roman{ize})}{\usecounter{ize}\setlength{\rightmargin}{\leftmargin}}
\item An Ehresmann connection on  $\pi_{XY}: Y \rightarrow X $  
\item A $GL(n)\times GL(k)$-equivariant global section
of  $\eta: \lvy \rightarrow \tauvf$
\item  A global trivialization 
 of $\lvy$ over $\tauvf$,
 $$ \Lambda: \lvy \rightarrow \tauvf \times \r^{k\times n}:
 w \mapsto (\eta(w),\lambda(w))$$ 
\hspace*{0.4in} 
where 
$\lambda: \lvy\rightarrow\rkn$
is a $\ga $-equivariant function
\item A symmetry 
 breaking $ \lambda: \lvy \rightarrow \r^{k\times n}$
where $\lambda^{-1}(0) \simeq \tauvf$  
\item  A global section of $JY\rightarrow Y$
\end{list}
\end{thm}

\bigskip

\begin{description}

\item[\noindent{\bf Proof } $(i)\Leftrightarrow(v):$] 
See GIMMsy.${}^{\ref{GIMMsy}}$

\item[$(ii)\Rightarrow(iv):$] The proof is motivated by Ref.~\ref{Bl}.
Let $\sigma:\tauvf\rightarrow\lvy$ be a  $GL(n)\times GL(k)$-equivariant section.
Since any element of $\eta^{-1}(u)$ can be expressed 
uniquely as $\sigma(u)\cdot A$ for some $A\in \rkn$, 
we may
define \mbox{$\lambda:\lvy\rightarrow\rkn$} by $\lambda(\sigma(u)\cdot A) = A$.
   Let  $B\,\in\,\rkn$. Then
\begin{eqnarray*}
\lambda((\sigma(u)\cdot B )\cdot(N,K,A))
& = & \lambda\left(\sigma(u)\cdot(N,K,BN+A)\right)\\
& = &  \lambda\left(\sigma(u\cdot(N,K))\cdot(K^{-1}BN + K^{-1}A)\right)\\
& = & K^{-1}BN + K^{-1}A \\
& = & (N,K,A)^{-1}\cdot \lambda(\sigma(u)\cdot B)\, .
\end{eqnarray*}
By definition, $\lambda(\sigma(u)\cdot A)= 0 \iff A=0$,
so $\lambda^{-1}(0) \simeq \tauvf$.

\item[$(iii)\Rightarrow(ii):$] 
This proof also is motivated by Ref.~\ref{Bl}.
Define $\sigma:\tauvf\rightarrow\lvy$ by $\sigma(u) = \Lambda^{-1}(u,0)$.
Using~(\ref{eqn:tauvfrkn}),
$$
\Lambda\left(\sigma(u)\cdot(N,K,0)\right)
 =  (u,0)\cdot (N,K,0) \\
 =  (u\cdot(N,K),0) \\
 =  \Lambda(\sigma(u\cdot(N,K))) \, .
$$
Since $\Lambda$ is bijective, 
$\sigma(u\cdot(N,K)) = \sigma(u)\cdot(N,K,0)$.

\item[$(i)\Rightarrow(iv):$]
Let $\gamma_y: T_yY\rightarrow V(T_yY)$ be an Ehresmann connection.
Define a map 
$$\lambda_\gamma : \lvy\rightarrow\rkn: \ptlvy \mapsto 
\epsilon^A |_{V(T_yY)}(\gamma(e_i))E^i_A$$
where $\{E^i_A\}$ is the standard basis of $\rkn$.
 Now,
if $v\in V(T_yY)$ then $e^i(v) = 0$.  So, 
\begin{eqnarray*}
\lefteqn{\lambda_\gamma(\ptlvy\cdot(N,K,A)) }\hspace{0.75in}\\
& = & \left((K^{-1})^B_A\epsilon^A- (K^{-1}AN^{-1})^B_ke^k\right) |_{V(T_yY)}\left(
\gamma(e_jN^j_i + \epsilon_AA^A_i)\right)E^i_B\\
& = & \left((K^{-1})^B_A\epsilon^A |_{V(T_yY)}(\gamma(e_j))N^j_i + (K^{-1})^{B}_CA^C_i\right) E^i_B\\
& = & (N^{-1},K^{-1},-K^{-1}AN^{-1})\cdot(\epsilon^B|_{V_y}\gamma(e_i))E^i_B\\
& = & (N,K,A)^{-1}\cdot\lambda_\gamma\ptlvy \; .
\end{eqnarray*}

\item[$(iv)\Rightarrow(i):$] Let $\lambda\ptlvy =\lambda^B_j\ptlvy E^j_B$.
Define a linear operator $\gamma_y$ on $T_yY$ by 
  $\gamma_y(\epsilon_A) := \epsilon_A$ and $\gamma_y(e_j) :=\lambda^B_j\ptlvy\epsilon_B$. 
  Then
$\gamma_y$ is the projection onto $V(T_yY)$ along  $\mbox{span}\{e_j-\lambda^B_j\ptlvy\epsilon_B\}$.  
The definition of $\gamma$ is  independent of choice of a point in $\pi_{Y\,\lvy}^{-1}(y)$.  
Indeed, if 
$(y,\{e_i',\epsilon_A'\})\in \pi_{Y\,\lvy}^{-1}(y)$ then 
$(y,\{e_i',\epsilon_A'\}) 
= \ptlvy\cdot(N,K,A)
= (y,\{e_jN^j_i + \epsilon_B A^B_i,\epsilon_BK^B_A\})$
for some $(N,K,A) \in \ga$.  
Define $\gamma_y'$ by 
$\gamma_y'(\epsilon_A') := \epsilon_A'$
 and $\gamma_y'(e_j') :=\lambda^B_j(y,e_i',\epsilon_A')\epsilon_B'$.
Hence $\epsilon_BK^B_A = \gamma_y'(\epsilon_BK^B_A) = \gamma_y'(\epsilon_B)K^B_A$
and thus 
$\epsilon_A = \gamma_y'(\epsilon_A)$.
Also, since $\lambda$ is $\ga $-equivariant, 
\begin{eqnarray*}
\gamma_y'(e_k') 
& = & \lambda^B_k(\ptlvy\cdot(N,K,A))\epsilon_B \\
& = & \left((K^{-1})^D_C \lambda^C_j\ptlvy N^j_k + (K^{-1})^A_CA^C_k \right) \epsilon_BK^B_D \\
& = & \gamma_y(e_j)N^j_k + \epsilon_B A^B_k \\ 
& = & \gamma_y(e_j)N^j_k + \gamma_y'(e_k') - \gamma_y'(e_j)N^j_k \; ,
\end{eqnarray*}
so $\gamma_y'(e_i) = \gamma_y(e_i)$. 
Therefore $\gamma_y'$ and $\gamma_y$
coincide on the frame $(y,\{e_i,\epsilon_A\})$.

\item[$(iv)\Rightarrow(iii):$] Define $\Lambda(w) := (\eta(w),\lambda(w))$.
Using~(\ref{eqn:etaequiv}) and~(\ref{eqn:tauvfrkn}), 
\begin{eqnarray*}
\Lambda(w\cdot(N,K,A)) 
& = & \left(\eta((w\cdot AN^{-1})\cdot(N,K,0)), \lambda(w\cdot(N,K,A))\right)\\ 
& = & \left(\eta(w)\cdot(N,K), (N,K,A)^{-1}\cdot\lambda(w)\right)\\ 
& = & \Lambda(w)\cdot(N,K,A)\, . 
\end{eqnarray*}
The transitive $\ga $--action on $\rkn$ 
defined in~(\ref{eqn:gaaction}) ensures that  
$\lambda$ is surjective.
To show that $\Lambda$ is surjective, let $(\eta(w),B) \in \tauvf \times \rkn$. 
If $C\in \rkn$ then  
$\eta(w\cdot C) = \eta(w)$ 
and 
$\lambda(w\cdot C) 
 =  (I_n, I_k, C)^{-1}\cdot \lambda(w) =  \lambda(w) + C$.
If $C= B-\lambda(w)$ then
$\Lambda(w\cdot C) = (\eta(w),B)$.
To show that $\Lambda$ is injective,
let $\Lambda(w) = \Lambda(w')$
for some $w, w' \in\lvy$.
 Then $w'\in\, \eta^{-1}(\eta(w))$, so $w' = w\cdot A$
 and $\lambda(w) = \lambda(w\cdot A) = (I_n,I_k,A)^{-1}\cdot \lambda(w).$
Hence, $A=0$.
Finally, note that
$\Lambda$ is $\rkn$-equivariant since $\rkn \subset \ga $.
~\blob
\end{description}

\begin{cor} 
An Ehresmann connection on ${\pi_{XY}}: Y\rightarrow X$
induces a flat connection  on 
 $ \eta: \lvy \rightarrow \tauvf$. 
\end{cor}

\proof If $\lambda_\gamma$ denotes the $\ga $-equivariant function
 obtained from an Ehresmann connection $\gamma$ in the  proof of 
Theorem~\ref{thm:five} ($(i)\Rightarrow(iv)$),
then  $R^*_A d\lambda_\gamma = d(\lambda_\gamma + A) = d\lambda_\gamma
 = Ad(A^{-1})\cdot d \lambda_\gamma$
 for
$A\in\,\rkn$. 
Let 
  $A= A^B_iE^i_B\in\,\rkn$ and $w= \ptlvy\in\,\lvy$.
Let $A^*_w$ be the fundamental vertical vector${}^{\ref{KN}}$ for $A$.
Then,
\begin{eqnarray*}
d\lambda_\gamma (A_w^*) & = & d\lambda_\gamma \left( \left. \frac{d}{dt} \,
\ptlvy\cdot \exp tA \,\right|_{t = 0}\right) \\
 & = & \left.\frac{d}{dt} \,\lambda_\gamma (y,\{e_i + tA^B_i\epsilon_B,\epsilon_A\} ) \,\right|_{t=0} \\
 & = & \left.\frac{d}{dt}\,(\epsilon^C-tA^C_je^j)\gamma(e_i + tA^B_i\epsilon_B)\,\right|_{t=0} E^i_C \\
 & = & \left.\frac{d}{dt}\epsilon^B\gamma(e_i) + t A^B_i\,\right|_{t=0} E^i_B \\
&  = & A \, .
\end{eqnarray*}
The connection $d\lambda_\gamma$ is flat 
because the Lie algebra $(\rkn, +)$ is Abelian.~\blob 

\bigskip

The pullback of the  {\em $(n+k)$-symplectic structure} $d\theta$ 
via  inclusion $i:\lvy\hookrightarrow LY $
is closed and  nondegenerate on $\lvy$,
 just as $d\theta$ is on
$LY$. 
The 
{\em ${(n+k)}$-symplectic structure equation}  for $\lvy$ is 
\beq\label{eqn:nkstruc}
d\hat f = -X_{\hat f}\hook i^*d\theta\; ,
\eeq
where 
 $\hat f\in\, C^\infty(\lvy,\r^{n+k}) $ and
$X_{\hat f}\in{\cal X}(\lvy)$, the vector space of vector fields of $\lvy$. 
Like the $p=1$ case of equation~(\ref{eqn2:struc1}),
equation~(\ref{eqn:nkstruc}) admits neither all vector fields
nor all $\rnplusk$-valued functions.

\bigskip

\noindent
{\bf Definitions} \begin{description}
\item[{$\bullet$}] $T^1(\lvy)$ is the vector space of {\em tensorial} 
$\r^{n+k}$-valued functions on $\lvy$.
\item[{$\bullet$}] The Lie algebra of {\em locally Hamiltonian vector fields} on $\lvy$
 is   
$$LHV^1(\lvy) := \{X\in\,{\cal X}(\lvy) \;|\; {\cal L}_X i^*d\theta = 0 \}\, .$$
\item[{$\bullet$}]  The vector space of {\em allowable  Hamiltonian observables} on $\lvy$
is 
$$HF^1(\lvy)  :=
 \{\hat f\in\, C^\infty(\lvy,\r^{n+k}) \;|\; d\hat f = -X\hook i^*d\theta, \; X\,\in\,LHV^1(\lvy)\}\, .
$$
\item[{$\bullet$}] The vector space of  {\em   Hamiltonian vector fields} on $\lvy$
 is   
$$HV^1(\lvy) := \{ X\in\,{\cal X}(\lvy) \;|\; X\hook i^*d\theta =  -d\hat f ,\,
\hat f \in HF^1(\lvy)\}\, .$$
\end{description}

Local coordinates on $\lvy$ are  $\{x^i,\,y^A,\,\pi^i_j,\,\pi^A_B,\pi^A_i\}$,
where 
$ \pi^i_j := e^i(\basisx j)$, $ \pi^A_B := \epsilon^A(\basisy B)$ and 
$ \pi^A_j := \epsilon^A(\basisx j)$.
In these coordinates,
$$i^*d\theta = d\pi^i_j\wedge dx^j \otimes \hat r_i + 
(d\pi^A_i\wedge dx^i +  d\pi^A_B \wedge dy^B )\otimes\hat s_A \, .$$
An element of $T^1(\lvy)$ can be expressed in local coordinates as 
$$
\hat f = f^i(x^k,y^C)\pi^j_i\otimes\hat r_j 
+ f^A(x^k,y^C)\pi^B_A\otimes\hat s_B
+ f^i(x^k,y^C)\pi^B_i\otimes\hat s_B \, .
$$

For $\hat f \in T^1(\lvy)$ we solve 
equation~(\ref{eqn:nkstruc}) 
locally  for $X_{\hat f}$.   This yields
\beq \label{eqn:xflvy}
X_{\hat f}  =  f^i\basisx i + f^A\basisy A   
- \partialx {f^k} j \pi^i_k\basispi i j  -  \partialy {f^C} B \pi^A_C
\basispi A B
-\left( \partialx {f^j} i \pi^A_j + \partialx {f^B} i \pi^A_B  \right)
\basispi A i \, ,   
\eeq
subject to the  constraints on $\hat f$,
\beq\label{eqn:constraint}
\partialy f^i A  = 0 \quad
\forall i = 1,\dots,n \;\;\mbox{and}\;\; A= 1,\dots,k.
\eeq
Thus, $T^1(\lvy)\not\subset HF^1(\lvy)$. 
Define
$$T^1_V(\lvy) := HF^1(\lvy)\cap  T^1(\lvy)\; .$$
Since $\lvy \subset LY$,
a point $w\in\,\lvy$ is a linear isomorphism
$w:\r^{n+k}\rightarrow T_{\pi_{Y\,\lvy}(w)}Y\, .$
Moreover, for $g \in \ga , (w\cdot g)V = w (g V)$ for each $V\in \rnplusk$.
Hence  $\hat f \in T^1(\lvy)$ 
corresponds bijectively to
${\f}\in\,{\cal X}Y$ by the relation $\hat f(w) = w^{-1}({\f}
 (\pi_{Y\,\lvy} (w)))$.
Now, if $\hat f\in\,T^1_V(\lvy)$ then,  by~(\ref{eqn:constraint}),
$\hat f$ is induced from a projectable vector field. 
Conversely, ${\boldmath f} \in \projectable$ gives a tensorial function $\hat f$
satisfying~(\ref{eqn:constraint}).
Thus,  {\em
$\;T^1_V(\lvy)$ is in bijective correspondence with $ \projectable$}.

 For $n\geq 2$ and $k\geq 2$
a vector field $X$ on $\lvy$
satisfies ${\cal L}_X(i^*d\theta) = 0$
if and only if in local coordinates
\begin{eqnarray}
X &=&  g^i(x^k)\basisx i +  g^A(x^k,y^C)\basisy A - \left(g^D,_{B}(x^k,y^C)\pi^A_D
+ \zeta^A,_{B}(x^k,y^C)\right) \basispi A B  \nonumber \\ 
&& \mbox{}
 - \left(g^i,_{l}(x^k)\pi^j_i + \xi^j,_{l}(x^k)\right)\basispi j l \label{eqn:hvlvy}\\
&& \mbox{} 
 - \left(g^B,_{i}(x^k,y^C)\pi^A_B + g^j,_{i}(x^k)\pi^A_j + \zeta^A,_{i}(x^k,y^C)
 + \eta^A,_{i}(x^k)
\right) \basispi A i \, , \nonumber
\end{eqnarray}
where $g,_i = \partialx g i$ and $g,_B = \partialy g B$.
Observe that $\pi_{Y\,\lvy*}X \in \projectable$.
From (\ref{eqn:nkstruc}) and (\ref{eqn:hvlvy})
we obtain the local expression for $\hat g \in\, HF^1(\lvy)$.
\beq\label{eqn:hflvy}
\hat g  =  \left( g^j(x^k)\pi^i_j + \xi^i(x^k)\right)\otimes \hat r_i 
+ \left( g^i(x^k)\pi^B_i + g^A(x^k,y^C)\pi^B_A + \zeta^B(x^k,y^C)\right)\otimes \hat s_B 
\; .
\eeq
As a result of equation~(\ref{eqn:hflvy}),
if  $n \geq 2$ and $k \geq 2$ then
$$HF^1(\lvy) \simeq T^1_V(\lvy)\oplus C^\infty(X,\r^n)\oplus C^\infty(Y,\r^{k})\, .$$

\noindent Observe that the kernel of 
the re\-pre\-sen\-ta\-tion $HF^1(\lvy) \rightarrow HV^1(\lvy): \hat f \mapsto \xhat f$
is the set of 
 the constant $\rnplusk$-valued functions.  Thus,
$$
HV^1(\lvy) \simeq HF^1(\lvy) / \rnplusk
$$
in analogy to the result${}^{\ref{No1}}$ on the full frame bundle $LY$
and to equation~(\ref{eqn:HVtoHFZ}).

\section{\normalsize\bf  ASSOCIATING THE MULTIVELOCITY SPACES TO $\lvy$}
\label{sec:linmodel}
   \setcounter{equation}{0}

Define a linear left action of $\ga $ on $\rkn$  by 
\beq\label{eqn:linleftaction}
(N,K,A)\odot W := KWN^{-1}\; .
\eeq
(We'll reserve the more conventional ``$\;\cdot\;$'' for the nonlinear
left action in (\ref{eqn:gaaction}).)

Use~(\ref{eqn:galvy}) and (\ref{eqn:linleftaction})
to define an equivalence relation on $\lvy\times \rkn$,
\beq\label{eqn:equivlin}
(\ptlvy,W) \sim 
\left(\ptlvy \cdot (N,K,A)\, , \, (N,K,A)^{-1}\odot W\right) 
\; .
\eeq
Equivalence classes $\lb2 \ptlvy, W\rb2$ are points in the
associated bundle $\lvy\times_\ga \rkn$.
(We'll reserve the more conventional ``$[\; , \;]$'' for 
a different equivalence class.)

\begin{lemma}\label{lemma:hom} The map
$$\psi:\lvy\times_\ga\,\rkn\rightarrow\mbox{\em Hom}_Y(TX,V(TY))\;\, : \;\, 
\lb2 \ptlvy, W\rb2\mapsto (y, -W^B_j({\pi_{XY}}_*e)^j\otimes\epsilon_B)$$ is a
vector bundle
isomorphism over $Y$, where  $\{({\pi_{XY}}_*e)^j\}$ is the coframe of $X$
at $\pi_{XY}(y)$
 dual to the frame $\{(\pi_{XY*}e)_j\}:= \{\pi_{XY*}e_j\}$ and 
$W = W^B_jE^j_B \in \rkn$.
\end{lemma}

\proof  By~(\ref{eqn:equivlin}),
$\psi$ is well defined.  
Since $\hbox{Hom\,}_Y(TX,V(TY))$ is vector bundle isomorphic to
 $T^*X\otimes_Y V(TY)$,
we represent $-W_j^B  ({\pi_{XY}}_*e)^j\otimes\epsilon_B$ as an element of  
$\lin(T_{{\pi_{XY}}(y)}X,V(T_yY))$.
Let 
$\psi\lb2\ptlvy,W\rb2 $  $= $ 
$\psi\lb2(y',\{e_i',\epsilon_A'\}),W'\rb2\,$.
Then $y = y'$ and there exists an 
$(N,K,A) \in \ga$ such that
$ (y',\{e_i',\epsilon_A'\})$
$= \ptlvy\cdot (N,K,A)$.
Subsequently,
$\psi\lb2(y',\{e_i',\epsilon_A'\}),W'\rb2$
$= $
$ \psi\lb2\ptlvy,KW'N^{-1}\rb2\,$. 
So $W = KW'N^{-1}$,
 and thus $W'= (N,K,A)^{-1}\odot W$. 
 So  $((y',\{e_i',\epsilon_A'\}),W') \sim  (\ptlvy,W)$,
and thus $\psi$ is injective.
Now 
$\psi$
 is linear on each fiber over $y$
and dim$(\lin(T_{{\pi_{XY}}(y)}X,V(T_yY)))$ $ = kn $ $= \mbox{dim}\left((\lvy)_y\times_\ga\rkn\right)$.
Thus $\psi$ is surjective on each fiber over $y$.
~\blob

\bigskip

Use~(\ref{eqn:galvy}) and~(\ref{eqn:gaaction})
to define a different equivalence relation on $\lvy\times \rkn$, 
\beq\label{eqn:equivkn}
(\ptlvy,W) \sim 
\left(\ptlvy \cdot (N,K,A)\, , \, (N,K,A)^{-1}\cdot W\right) 
\, .
\eeq
Equivalence classes  $[\ptlvy,W]$ are points in the associated bundle
$(\lvy\times_\ga \rkn)_{\mbox{Aff}}\,$.
Ragionieri and Ricci${}^{\ref{RR}}$ show that $JY$ is an affine bundle
over $Y$ modeled on $\hom_Y(TX,V(TY))$.  
We now can reproduce this result with bundles associated to $\lvy$.

\begin{thm}
The bundle  $(\lvy\times_\ga \rkn)_{\mbox{\em Aff}}\,$ is an affine
bundle over $Y$ with underlying vector bundle $\lvy\times_\ga\rkn$,
and  the map
$$
\begin{array}{rcl}
\hat\psi :(\lvy\times_\ga \rkn)_{\mbox{Aff}} & \rightarrow  &  JY  \\ 
 \mbox{$[$}\ptlvy , W ] &\mapsto &\left( y,-W^B_j ({\pi_{XY}}_*e)^j\otimes\epsilon_B + ({\pi_{XY}}_*e)^j\otimes e_j\right)
\end{array}
$$ 
is an affine bundle
isomorphism onto the jet bundle $JY$.
\end{thm}

\proof Define a difference function $\hat\delta$ on 
$(\lvy\times_\ga \rkn)_{\mbox{Aff}}\,$
on each fiber over $Y$ by 
$$ \hat\delta_y([w,W_1]\, , \, [w,W_2]) = \lb2 w\, , \,W_1 -  W_2\rb2$$
where $y = \pi_{Y\,\lvy}(w)$.  
The range is $(\lvy)_y\times_\ga\rkn$.
The function $\hat\delta_y$ is well defined.  Indeed, if
 $w' \in \pi_{Y\,\lvy}^{-1}(\pi_{Y\,\lvy}(w))$, 
then there exists
$(N,K,A)\in \ga $ such that $w' = w \cdot (N,K,A)$.
Let $[w', W_1'] = [w,W_1]$ and $[w', W_2'] = [w,W_2]$.
   Then $W_1 = (N,K,A)\cdot W_1'$
and $W_2 = (N,K,A)\cdot W_2'$ and thus
\beqann
\lb2 w' , W_1' - W_2'\rb2
 & = &  \lb2 w \cdot (N,K,A) ,
   K^{-1} W_1 N + K^{-1}A  -   (   K^{-1} W_2 N + K^{-1}A)\rb2 \\
& = & \lb2 w\cdot(N,K,A), (N,K,A)^{-1}\odot (W_1-W_2)\rb2 \\
& = & \lb2  w , W_1 - W_2\rb2 \; .
\eeqann
Similarly, we may verify the other properties of the difference
function.

Let $(y,\gamma_y) = \hat\psi[\ptlvy,W]$. 
Then 
 $\gamma_y \in\,\hbox{Lin}(T_{{\pi_{XY}}(y)}X,T_yY)$
and  
${\pi_{XY}}_*\circ \gamma_y 
 =  \;\mbox{Id}_{T_xX}$.
Therefore, $(y,\gamma_y) \in JY$.
Using~(\ref{eqn:equivkn}),
 we can
show that $\hat\psi$ is well defined.  
To show that $\hat\psi$ is surjective
choose an arbitrary point $\ptlvy \in \lvy$ and define $\gamma^B_j(y)$ by
$\gamma_y({\pi_{XY}}_*e_j) = e_j + \gamma^B_j(y)\epsilon_B$. Then 
 $(y,\gamma_y) = \hat\psi[\ptlvy,-\gamma^B_j(y)E^j_B]$.
To show that $\hat\psi$ is injective
let $\hat\psi[\ptlvy,W] = \hat\psi[(y',\{e_i',\epsilon_A'\}),W']$.
Then $y = y'$ and there is a unique $(N,K,A)\in\,\ga $ such that 
$(y',\{e_i',\epsilon_A'\}) = (y,\{e_i, \epsilon_A\})\cdot(N,K,A)$.
So 
$$
\hat\psi[(y,\{e_i,\epsilon_A\})\cdot(N,K,A)\, ,\,W'] 
 =  \hat\psi[(y,\{e_i,\epsilon_A\})\, ,\, (N,K,A)\cdot W'] 
\qquad\qquad\qquad\mbox{} 
$$
$$ \mbox{}\qquad\qquad\qquad\qquad
= \left(y, - ( KW'N^{-1}-AN^{-1})^B_j({\pi_{XY}}_* e)^j \otimes \epsilon_B 
+ ({\pi_{XY}}_*e)^j\otimes e_j \right) ,
$$
giving us 
$KW'N^{-1}-AN^{-1} = W$ or $W' = (N,K,A)^{-1}\cdot W$.
Thus $(\ptlvy,W) \sim ((y',\{e_i',\epsilon_A'\}),W')$.
Finally, 
let $\delta$
be the difference function on $JY$.
 Let $w= \ptlvy\in \lvy$, $\pi_{Y\,\lvy}(w) = y$
and let $W$ and $W' \in \rkn$. 
  Then,
$$
 \psi\left(\hat\delta([w,W]\, ,\, [w,W'])\right) 
 = \delta\left(\hat\psi[w,W]\, ,\, \hat\psi[w,W']\right) 
$$
where $\psi$ is defined in Lemma~\ref{lemma:hom}.~\blob

\section{\normalsize\bf ASSOCIATING THE MULTIPHASE SPACES TO $\lvy$}
\label{sec:affmodel}
   \setcounter{equation}{0}

Define a linear left action of $\ga$ on $\r^{n\times k}$ by
\begin{equation}\label{eqn:g-action}
(N,K,A)\odot B := NBK^{-1}\, ,
\eeq  
and define an equivalence relation on $\lvy\times\rnk$ similar 
to (\ref{eqn:equivlin}),
using  (\ref{eqn:galvy}) and (\ref{eqn:g-action}).
The resulting equivalence classes, denoted by $[ \ptlvy, B]_{Gu}$, are points in the associated bundle
$(\lvy\times_\ga\rnk)_{Gu}$. 
The following result comes
from modifying the proof of 
Lemma~\ref{lemma:hom}.

\begin{lemma}\label{lemma:j*y} Let  $\{E^C_j\}$ be the standard basis of $\rnk$
and  $B = B^j_CE^C_i \in \rnk$. Then  
\beqann
\rho_{Gu} : (\lvy\times_\ga \rnk)_{Gu} 
	& \rightarrow & \mbox{Hom\,}_Y(V(TY),TX)\\
\mbox{}  [\ptlvy,B]_{Gu}
	&\mapsto& (y,B^j_C\epsilon^C |_{V(T_yY)}\otimes f_j)
\eeqann
is a vector bundle isomorphism over $Y$.
\end{lemma}

\noindent If, instead of the action
defined in~(\ref{eqn:g-action}),
we use a {\em $GL(n)$-normalized} action
\beq\label{eqn:gnormed}
(N,K,A)\cdot B =\hbox{det}(N^{-1})NBK^{-1} \, ,
\eeq
then the bundle with fiber $\rnk$ associated to $\lvy$ 
is $\kt$.
Using~(\ref{eqn:galvy}) and~(\ref{eqn:gnormed}),
the new equivalence classes denoted by
$[\ptlvy,B]$ are points in the associated bundle $\lvy\times_\ga\rnk$.
Define 
\beqann
\rho_{KT} : \lvy \times_\ga \r^{n\times k}  
	&\rightarrow& \hom_Y(V(TY),\wedge^{n-1}X)  \\
\mbox{} [\ptlvy ,B] 
	&\mapsto &(y,B^j_C\epsilon^C|_{V(T_yY)}\otimes \omega({\pi_{XY}}_*e)_j)\, , 
\eeqann
where if $f$ is a set of $n$ linearly independent vectors  $\{f_i\}$ then
\beq\label{eqn:omega(f)}
\omega(f)  :=  f^1\wedge f^2\wedge \dots\wedge f^n \quad\mbox{and}\quad 
\omega(f)_j  :=  f_j\, \hook \omega(f)\;.
\eeq
 Like $\rho_{Gu}$ in Lemma~\ref{lemma:j*y},
 $\rho_{KT}$ is a vector bundle isomorphism over $Y$.

Define a linear left action of $\ga $ on the vector space $\r^{n\times k}\times \r$ as follows.
\beq\label{eqn:ga-action}
(N,K,A)\cdot (B,\lambda) := \hbox{det}(N^{-1})\left(NBK^{-1}, \lambda - \tr(BK^{-1}A)\right)
\eeq
The associated vector bundle
$\lvy\times_{\ga } (\r^{n\times k}\times \r)$ 
is constructed  using (\ref{eqn:galvy}) and (\ref{eqn:ga-action}).
Equivalence classes (points in the vector bundle) are denoted by
$[\ptlvy,(B,\lambda)]$.

\begin{thm}\label{thm:z}
The map 
$$\begin{array}{rcl}
\rho_{Z} : \lvy\times_{\ga } (\r^{n\times k}\times \r) &
\rightarrow & \wedge^n Y \\
\left[ \ptlvy,(B,\lambda)\right] & \mapsto & (y, B^j_B\epsilon^B
\wedge\omega(e)_j + \lambda\omega(e))
\end{array}
$$
is a vector bundle monomorphism over $Y$,
and the range of $\rho_{Z}$
is $Z$.
\end{thm}

\proof We may use (\ref{eqn:galvy}), (\ref{eqn:coframeaction}),
 (\ref{eqn:omega(f)}) and (\ref{eqn:ga-action}) 
to show that $\rho_{Z}$ is well defined.
Observe that $\omega(e)$ kills vertical vectors, and   
if $n \geq 2$ then
$\epsilon^A\wedge \omega(e)_i$ kills two
vertical vectors. 
 Thus, the range of  $\rho_{Z}$ is in $Z$.  
For any $z \in Z$, in local coordinates,
$$
z = \rho_{Z} \left[ \left(\pi_{YZ}(z), \{ \basisx i,
 \basisy A \}\right),\left((p^i_A(z)),p(z)\right)\right]
$$
where
 $p^i_A(z)$ and
$p(z)$ are the coordinates defined in~(\ref{eqn:zcoords}).
Also, $\rho_{Z}$ is linear on fibers over $Y$, $\rho_Z$ preserves fibers,
and the fibers have the same dimension.~\blob

\bigskip

Use $\Pr_1:\rnk \times \r$ to $\rnk$ and the maps
$\rho_{KT}$ and $\rho_{Z}^{-1} |_Z$ to define the projection
$\pi_{KT\,Z} : Z\rightarrow 
\hom_Y(V(TY),\wedge^{n-1}X)$ .
See diagram~(\ref{diag:assoc}).
\beq\label{diag:assoc}
\begin{array}{ccccc}
\lvy\times_{\ga}(\rnkxr)&{\rho_{Z}^{-1}|_Z\atop\longleftarrow}&Z\\
&&\\
\downarrow&&\downarrow\pi_{KT\,Z}\\\
&&\\
\lvy\times_\ga\rnk&
{\rho_{KT} \atop\longrightarrow}&\kt
\end{array}
\eeq

\noindent Using $\rho_{Z}$, the relationships between local canonical
 coordinates $\{x^i,y^A,p^j_B,p\}$ on $Z$
and local coordinates $\{x^i, y^A, \pi^j_k, \pi^B_C, \pi^D_l\}$ on
 $\lvy$ are
\beq\label{eqn:coordslvyz}
p^j_B = \det(\pi^l_m)\,B^i_A\pi^A_B(\pi^{-1})^j_i \; \quad\mbox{and}
\quad p = \det(\pi^l_m)(B^i_A\pi^A_k(\pi^{-1})^k_i + \lambda)\; .
\eeq
The first equation in~(\ref{eqn:coordslvyz})
also relates coordinates in $\hom_Y(V(TY),\wedge^{n-1}X)$ to 
those in $\lvy$.
 
The identification of the linear and affine multiphase spaces
with bundles associated to $\lvy$ allows us to recast certain known 
relationships between these spaces in our more general context.

\begin{thm}
The bundle $Z$ is an affine bundle over
$\kt$. (See Ref.~\ref{CCI}.) The vector bundle underlying $Z$ is 
$\kt \times_Y \pi_{YZ}^*\wedge^nX$. 
\end{thm}

\proof Define a difference function
$\delta : Z \times Z \rightarrow \kt \times_Y \pi_{YZ}^*\wedge^nX$
by $\delta ([u,(B,\lambda_1)],[u,(B,\lambda_2)]) = ([u,B],(\lambda_1-\lambda_2)\,\omega(u))$,
where if $u = \ptlvy$ then $\omega(u) = \omega(e)$.
The map $\delta$ is well defined
and preserves fibers over $\kt$.   To check that
$\delta$ is well defined,
let $(u',(B',\lambda'))\sim(u,(B,\lambda))$.  Then $u' = u\cdot g$ and 
$(B',\lambda_i') = g^{-1}\cdot (B,\lambda_i)$ for
some $g = (N,K,A)\in \ga$ and $i = 1,2$. 
Now, $\omega(u\cdot g) = \det(N^{-1})\,\omega(u)$, so,
\beqann
\delta([u', (B',\lambda_1')] \, , \, [u',(B',\lambda_2')])
 & = &
\delta([u\cdot g, g^{-1}\cdot(B,\lambda_1)] 
\, , \, [u\cdot g, g^{-1}\cdot(B,\lambda_2)]) \\
& = & \delta\left([u\cdot g, \det N(N^{-1}BK,\lambda_1 + \tr(BAN^{-1}))]\, , 
\right.\\
&&  \hspace*{0.8in}\,
 \left.[u\cdot g,\det N (N^{-1}BK,\lambda_2 + \tr(BAN^{-1}))]\right) \\
& = &  \left([u\cdot g, g^{-1}\cdot B]\, ,\,\det N (\lambda_1-\lambda_2)
\omega(u\cdot g)
\right) \\
& = &  \left([u,B]\, , \,  (\lambda_1-\lambda_2)\,\omega(u)\right) \\
& = & \delta ([u,(B,\lambda_1)]\, , \, [u,(B,\lambda_2)])\; . 
\eeqann
It is straightforward to verify the other properties of the difference function
$\delta$.~\blob

\begin{lemma}\label{lemma:secform}
Any section 
of the bundle $\pi_{KT\,Z}:Z \rightarrow \kt$
 may be expressed as a map
$[u,B] \mapsto [u,(B,\eta(u,B))]$, 
where $\eta: \lvy \times \rnk \rightarrow \r$ satisfies the equation
\beq\label{eqn:zseccond}
(\det N) \left(\eta(u,B) + \tr (BAN^{-1})\right) = \eta(u\cdot (N,K,A),(\det N)N^{-1}BK) 
\eeq
 for all $(N,K,A) \in \ga$.
Conversely, any map 
$\eta : \lvy \times \rnk \rightarrow \r$
that satisfies (\ref{eqn:zseccond}) induces a section
$[u,B] \rightarrow [u,(B, \eta(u,B))]$
of the bundle $\pi_{KT\,Z}:Z \rightarrow \kt$.

\end{lemma}

\proof By definition, a section $s: \kt \rightarrow Z$
must satisfy 
$s [u,B] = [u,(B,\lambda)]$ 
where $u \in \lvy$,  $B \in \rnk$ and
$\lambda \in \r$.  
For each $u \in \lvy$,
define $\sigma_u: Z \rightarrow \rnkxr:
[u,(B,\lambda)] \mapsto (B,\lambda)$.  
Because $Z\simeq \lvy\times (\rnk \times \r)$,
 $\sigma_u$ is well defined. 
Define
$\eta: \lvy \times \rnk \rightarrow \r $ by
$ \eta(u,B) := \Pr_2(\sigma_u(s[u,B]))$,
where $\Pr_2$ is projection in the second slot.  Then
$s[u,B] = [u,(B,\eta(u,B))]$.
Now
$s[u,B] = s[u\cdot g, g^{-1}\cdot B]$
for $g \in\ga$,
so $[u\cdot g, (g^{-1}\cdot B, \eta(u\cdot g,
g^{-1}\cdot B))] =
 [u,(B, \eta(u,B))] = [u\cdot g, g^{-1}\cdot(B,\eta(u,B))]$.
So $g^{-1}\cdot (B,\eta(u,B)) = (g^{-1}\cdot B, \eta(u\cdot g, g^{-1}\cdot B))$
and thus  we obtain equation~(\ref{eqn:zseccond}).
Conversely, if  $s$
 satisfies~(\ref{eqn:zseccond}) then 
 $s$ is well defined.~\blob

\begin{thm}(Sardanashvily$\,{}^{\ref{Sa}}$)\label{prop:split}
An Ehresmann connection on $Y$
is equivalent to a splitting of the short exact sequence of vector
bundle homomorphisms over $\mbox{Id}_{\,Y}$,
$$0\rightarrow {\pi_{XY}}^*(\wedge^n X)\rightarrow Z\rightarrow \kt\rightarrow 0\, .$$
\end{thm}

\proof Assume an Ehresmann connection on $Y$.
By Theorem~\ref{thm:five},
a connection is equivalent to a
$\ga$-equivariant map 
$\lambda: \lvy \rightarrow \rkn$.
Define a map
$$\eta:\lvy\times\rnk \rightarrow \r : (u,B) \mapsto \tr(B\lambda(u)) \, .$$
 Let $g=(N,K,A) \in \ga$.  Using~(\ref{eqn:gaaction}),
\beqann
(\det N) \left(\eta(u,B) + \tr(BAN^{-1})\right) 
& = &  (\det N)\, \tr\left(N^{-1}BK(K^{-1}\lambda(u)N +K^{-1}A )\right)\\
& = & \tr\left((\det N)N^{-1}BK\lambda(u\cdot g) \right) \\
& = &  \eta(u\cdot g, (\det N) N^{-1}BK)
\eeqann
so  $\eta$ satisfies 
(\ref{eqn:zseccond}). 
By Lemma~\ref{lemma:secform}, we may define a
section
$s : \kt \rightarrow Z$ and it is routine to check that $s$ is linear
on fibers over $Y$. 

Conversely, let the section $s: \kt \rightarrow Z$ be 
 linear on fibers over $Y$.
By Lemma~\ref{lemma:secform}, 
we write $s[u,B] = [u,(B,\xi(u,B))]$, and $\xi$ satisfies  
 $\xi(u,B_1+cB_2) = \xi(u,B_1)+ c\,\xi(u,B_2)$. 
For $u\in \lvy$ we define $\xi_u :\rnk \rightarrow \r$ 
by $\xi_u(B) := \xi(u,B)$.  Since $\xi_u$ is linear we may uniquely represent 
$\xi_u$ as $\xi_u(B) = \tr(BW(\xi_u))$ where $W(\xi_u) \in \rkn$. 
We now argue that $u \mapsto W(\xi_u)$ is 
a symmetry-breaking map and thus is equivalent to an
Ehresmann connection on $Y$.
 Indeed, the action of  $\ga$
 on $\rkn$ given in (\ref{eqn:gaaction}) is transitive,
so it remains only to show equivariance.
By equation~(\ref{eqn:zseccond}),
\beqann
\det N (\xi_u(B) + \tr (BAN^{-1})) &=& \xi_{u\cdot g}((\det N)N^{-1}BK) \; \; \forall \,B \;\\
\tr( BW(\xi_u))  + \tr(BAN^{-1}) 
& = & \tr( BKW(\xi_{u\cdot g})N^{-1}) \;\; \forall \, B\\
\tr\left(B( W(\xi_u) + AN^{-1} - KW(\xi_{u\cdot g})N^{-1})\right)
& = & 0  \;\;\forall B\\
KW(\xi_{u\cdot g})N^{-1} &=& W(\xi_u) + AN^{-1} \\
W(\xi_{u\cdot g}) & = & K^{-1}W(\xi_u)N + K^{-1}A \\
W(\xi_{u\cdot g}) & = & g^{-1}\cdot W(\xi_u) \; \; \blob
\eeqann
If $f: Z \rightarrow {\pi_{XY}}^*(\wedge^nX)$ splits 
the sequence in Theorem~\ref{prop:split}, then 
we may write 
$f[u,(B,\lambda)]$ $=[u,(0,\lambda-\xi(u,B))]$,
where $\xi$ is the map in the proof of  Theorem~\ref{prop:split}.
Let $\gamma$ be the connection equivalent to $f$.
 In local coordinates,
$$
f:\; p^i_Ady^A\wedge d^n x_i + p\, d^n x \mapsto (p + p^i_A\gamma^A_i)d^n x\; .
$$

\section{\normalsize\bf  AFFINE MULTISYMPLECTIC GEOMETRY FROM (n+k)-SYM\-PLEC\-TIC GEOMETRY}
\label{sec:affgeom}

 \setcounter{equation}{0}

We will construct 
the mul\-ti\-sym\-plec\-tic potential $\Theta$ on $Z$
from the $(n+k)$-symplectic potential 
$i^*\theta$ on $\lvy$.
First, some preliminary remarks about tensor-valued differential forms on $LY$
are necessary.

Let $\{R_\mu\}, \, \mu = 1,\dots,n+k$, be the standard basis of $\rnplusk$
and let $\{R^\mu\}$ be the corresponding dual basis.
For convenience, define $R_{\mu_1\cdots\mu_m} := R_{\mu_1}\wedge\cdots\wedge
R_{\mu_m} \in \wedge^m\rnplusk$ and 
$R^{\mu_1\cdots\mu_m} := R^{\mu_1}\wedge\cdots\wedge
R^{\mu_m} \in \wedge^m{\rnplusk}^*$.
 Let $\alpha$ be a $\wedge^r\rnplusk$-valued $p$-form and let $\beta$
be a $\wedge^s\rnplusk$-valued $q$-form on a manifold.  
Then 
$\alpha = \alpha^{\mu_1 \cdots \mu_r}\otimes R_{{\mu_1 \cdots \mu_r}}$ 
and
$\beta = \beta^{\mu_1 \cdots \mu_s}\otimes R_{{\mu_1 \cdots \mu_s}}$.
Define 
$$
\alpha\wedge\beta := (\alpha^{\mu_1 \cdots \mu_r}\wedge\beta^{\nu_1 \cdots \nu_s}
)\otimes R_{\mu_1 \cdots \mu_r\nu_1 \cdots \nu_s} \, .
$$
Observe that
$\alpha\wedge\beta = (-1)^{pq + rs}\beta\wedge\alpha \, .$
Let $\theta$ be the canonical soldering form on $LY$.
Define the  $\wedge^m\rnplusk$-valued $m$-form $\wedge^m\theta$ on $LY$ by
$\wedge^0\theta := 1$ and 
$$
\wedge^m\theta := \overbrace{\theta\wedge\cdots\wedge\theta}^{m} \quad
 \mbox{if} \quad m\geq 1\; . 
$$
Note that
\beq\label{eqn:dmtheta}
d(\wedge^m\theta)  =  m \, d\theta\wedge(\wedge^{m-1}\theta) \, .
\eeq

\begin{lemma}\label{lemma:mtheta}
If $1 \leq m \leq n+k-1$ then  $d(\wedge^m\theta)$ is closed and nondegenerate.
\end{lemma}

\proof If $m=1$ then the above form is $d\theta$.
Let $2 \leq m \leq n + k-1$.
The form 
 is closed because $d^{\,2}(\wedge^m\theta) = 0$.
Let $X \in {\cal X}(LY)$ satisfy $0 = X\hook d(\wtheta m)$.
By equation~(\ref{eqn:dmtheta}),
 $0 = X\hook (d\theta\wedge(\wtheta {m-1}))$.
We may introduce a torsionless connection $\omega$ on $LY$
without loss of generality.  Thus $0 = D\theta = d\theta + \omega\wedge\theta$,
or $$d\theta^\mu = -\omega^\mu_\nu\wedge \theta^\nu\, .$$ 
At $w \in \lvy$,  we can write 
 $X_w = X^\mu B_{\mu}(w) + X^\mu_\nu E^{*\mu}_{\nu}(w)$,
where $\{E^{*\mu}_\nu (w)\}$ is the frame of fundamental vertical vector fields
and 
$\{ B_\mu (w)\}$ is  the horizontal frame complementary to  $\{ E^{*\mu}_\nu (w)\}$.
So, 
\beqa
0 & = & X_w\, \hook (\theta^\nu\wedge\omega^{[\mu_1}_\nu\wedge\theta^{\mu_2}\wedge
\cdots\wedge\theta^{\mu_m]}) \nonumber \\
0  & = & X^\nu \omega^{[\mu_1}_\nu\wedge\theta^{\mu_2}\wedge
\cdots\wedge\theta^{\mu_m]} 
-\theta^\nu\wedge X^{[\mu_1}_\nu \theta^{\mu_2}\wedge
\cdots\wedge\theta^{\mu_m]} 
+ \theta^\nu \wedge \omega^{[\mu_1}_\nu
  X^{\mu_2}\theta^{\mu_3}\wedge\cdots\wedge \theta^{\mu_m]}
\nonumber \\
&  & \qquad + \sum^m_{i=3} (-1)^i \theta^\nu \wedge \omega^{[\mu_1}_\nu
  \wedge \theta^{\mu_2}\wedge\cdots\wedge
   X^{\mu_i}\widehat{\theta^{\mu_i}}\wedge\cdots\wedge \theta^{\mu_m]} \nonumber \\
0  & = & X^\nu \omega^{[\mu_1}_\nu\wedge\theta^{\mu_2}\wedge
\cdots\wedge\theta^{\mu_m]} 
-\theta^\nu\wedge X^{[\mu_1}_\nu \theta^{\mu_2}\wedge
\cdots\wedge\theta^{\mu_m]} \nonumber \\
&  & \qquad + \sum^m_{i=2} (-1)^i \theta^\nu \wedge 
X^{[\mu_i}\omega^{\mu_1}_\nu
  \wedge \theta^{\mu_2}\wedge\cdots\wedge\widehat{\theta^{\mu_i}} \wedge
\cdots\wedge \theta^{\mu_m]} \nonumber 
\eeqa
where $\widehat{\theta^{\mu_i}}$ denotes the omission of $\theta^{\mu_i}$.
Now use $[\mu_i,\mu_1,\dots,\widehat{\mu_i},\dots,\mu_m] =(-1)^{i-1}[\mu_1,\dots ,\mu_m]$ to combine terms, obtaining
\beqa
0 & = & X^\nu \omega^{[\mu_1}_\nu\wedge\theta^{\mu_2}\wedge
\cdots\wedge\theta^{\mu_m]} 
-\theta^\nu\wedge X^{[\mu_1}_\nu \theta^{\mu_2}\wedge
\cdots\wedge\theta^{\mu_m]} \nonumber \\
& & \mbox{}\qquad
- (m-1) \theta^\nu\wedge X^{[\mu_1}
\omega^{\mu_2}_\nu\wedge\theta^{\mu_3}\wedge
\cdots\wedge\theta^{\mu_m]} \, .  \label{eqn:(*)}
\eeqa
 Now  evaluate both sides of equation~(\ref{eqn:(*)}) 
at 
$(B_{\alpha_1},\dots,B_{\alpha_m})(w)$ to obtain
\beqann
0 & = & -\theta^\nu \wedge X_\nu^{[\mu_1}\theta^{\mu_2}\wedge \cdots \wedge
	\theta^{\mu_m]}(B_{\alpha_1},\dots,B_{\alpha_m})(w) \\
0 & = & (X_\nu^{[\mu_1}\theta^{|\nu|}) \wedge\theta^{\mu_2}\wedge \cdots \wedge
	\theta^{\mu_m]}(B_{\alpha_1},\dots,B_{\alpha_m})(w) \\
0 & = & \sum_{\sigma \in \mbox{\tiny Perm}(m)}\mbox{sgn}(\sigma)\,
	\theta^{\nu} \otimes    
X_{\nu}^{[\mu_{\sigma(1)}}\theta^{\mu_{\sigma(2)}}\otimes \cdots \otimes
	\theta^{\mu_{\sigma(m)}]}(B_{\alpha_1},\dots,B_{\alpha_m})(w) \\
0 & = & \sum_{\sigma\in \mbox{\tiny Perm}(m)}\mbox{sgn}(\sigma)\, 
	X_{\alpha_1}^{[\mu_{\sigma(1)}} \delta^{\mu_{\sigma(2)}}_{\alpha_2}\cdots 	
	\delta^{\mu_{\sigma(m)}]}_{\alpha_m} \\
0 & = & X^{[\mu_1}_{\alpha_1} \delta^{\mu_2}_{\alpha_2}\cdots
	\delta^{\mu_m]}_{\alpha_m} \, .
\eeqann
 Choose a sequence $\{\mu_I\}$ such that the
$\mu_I$ are distinct for $2 \leq I \leq m$ and 
choose a  sequence $\{\alpha_I\}$  such that 
$ \alpha_I = \mu_I $ for  $2\leq I \leq m$.  Thus $X^{\mu_1}_{\alpha_1} =0$,
for any choice of $\mu_1$ and $\alpha_1$.
Return to~(\ref{eqn:(*)}) and evaluate
 both sides of the equation at
$(E^{*\alpha}_\beta, B_{\alpha_2}, \dots, B_{\alpha_{m}})$ 
Thus, by a calculation similar to that above,
\beq\label{eqn:deltas}
0  =  X^\alpha \delta^{[\mu_1}_{\beta} \delta^{\mu_2}_{\alpha_2}
\cdots \delta^{\mu_m]}_{\alpha_m}
  + (m-1)\delta^\alpha_{\alpha_2}X^{[\mu_1}
 \delta^{\mu_2}_{\beta}\delta^{\mu_3}_{\alpha_3} \cdots \delta^{\mu_m]}_{\alpha_m}
\, .
\eeq
Again, choose a
sequence $\{\mu_I\}$ such that the $\mu_I$ are distinct if $2\leq I \leq m$. 
Now let $\alpha_I = \mu_I$ for $2 \leq I \leq m$ and let 
$\alpha= \alpha_2 =\mu_2 =\beta = \mu_1$.
Then the second term in equation~(\ref{eqn:deltas}) vanishes, else 
$\alpha= \alpha_2 = \beta$ which contradicts the hypothesis that the $\mu_i$ are distinct.~\blob

\bigskip

Pulling back $d(\wedge^m\theta)$ via inclusion 
$i: \lvy\hookrightarrow LY$,
it follows that 
 $d(\iwtheta m)$ is a closed, nondegenerate, $\wedge^m\rnplusk$-valued form on $\lvy$ if $1 \leq m\leq n+k-1$.
When expressing the standard basis of $\rnplusk$
as $\{\hat r_i, \hat s_A\}, \, i = 1,\dots,n,\, A= 1,\dots,k$
instead of $\{R_\mu\}, \mu = 1,\dots,n+k$,
but use upper case and lower case Latin indices in place
of lower case  Greek indices.
Then
 $i^*\theta = \theta^i\otimes \hat r_i + \theta^A\otimes \hat s_A$,
and, so by induction on $m$, 
 \beq\label{eqn:binomial}
\iwtheta m = \sum^m_{l=0}\left({m\atop l}\right)
\theta^{A_1}\wedge\cdots\wedge\theta^{A_l}\wedge\theta^{i_1}\wedge\cdots\wedge
\theta^{i_{m-l}}R_{A_1\cdots A_l i_1\cdots i_{m-l}}\, .
\eeq

Let $m=n$. If
$\sigma$ is a permutation of $n$ elements expressed by 
$(\sigma(1),\dots ,\sigma(n)) =(i_1,\dots ,i_n )$
then let $\epsilon_{i_1\cdots i_n}$ denote the sign of $\sigma$.
Let $\ptlvy \in \lvy$.
Using~(\ref{eqn:omega(f)}),
\beq\label{eqn:omega}
\omega(e) = \frac{1}{n!}\epsilon_{i_1\cdots i_n}e^{i_1}\wedge\cdots\wedge e^{i_n}
\eeq
and 
\beq\label{eqn:omegaj}
\omega(e)_j = \frac{1}{(n-1)!}\epsilon_{ji_1\cdots i_{n-1}}
e^{i_1}\wedge\cdots\wedge e^{i_{n-1}}.
\eeq
Define
$$
\phibl :  \lvy \rightarrow  Z :
 w \mapsto \rho_{Z}[w,(B,\lambda)]
$$
 The map $\phibl$ is fiber preserving over $Y$.

\begin{thm}\label{thm:nthetaV} Let $n\geq 2$.
The $\wedge^n\rnplusk$-valued $n$-form $\wedge^n i^*\theta$ on $\lvy$
can be related to the canonical $n$-form $\Theta$ on $Z$ by 
\beq\label{eqn:n-forms} 
\left< \wedge^n i^*\theta , \mbox{\boldmath $V$}(B,\lambda) \right> = \phi_{(B,\lambda)}^*\Theta
\eeq
where the map $\mbox{\boldmath $V$} : \rnk\times \r \rightarrow  \wedge^n{\r^{n+k}}^*$ has components 
\beqann
V_{i_1\dots i_n}(B,\lambda) &=&\frac{1}{n!} \lambda\epsilon_{i_1\dots i_n}\, ,\\
 V_{Ai_1\dots i_{n-1}}(B,\lambda)&=& \frac{1}{n!}B^j_{A}\epsilon_{ji_1\dots i_{n-1}}\\
 \mbox{and} \qquad V_{A_1\dots A_li_1\dots i_{n-l}}(B,\lambda)
      &=& 0 \;\; \forall \, l \geq 2. 
\eeqann
\end{thm}

\proof
Using equations~(\ref{eqn:binomial}), (\ref{eqn:omega}) and (\ref{eqn:omegaj}),
 for $w = \ptlvy \in \lvy$,
\beqann
\lefteqn{
\left< \wedge^n i^*\theta , \mbox{\boldmath $V$}(B,\lambda) \right> (w) }
\hspace{.75in}\\
& = & \frac{1}{n!}\lambda\epsilon_{i_1\cdots i_n} \theta^{i_1}\wedge\cdots\wedge
\theta^{i_n}(w) + \frac{n}{n!}B^j_A\epsilon_{j{i_1\cdots i_{n-1}}}  
\theta^A\wedge\theta^{i_1}\wedge\cdots\wedge\theta^{i_{n-1}}(w) \\
& = & \left(\frac{1}{n!}\lambda\epsilon_{i_1\cdots i_n}e^{i_1}\wedge\cdots\wedge e^{i_n} + \frac{1}{(n-1)!}B^j_A\epsilon_{j{i_1\cdots i_{n-1}}}  
\epsilon^A\wedge e^{i_1}\wedge\cdots\wedge e^{i_{n-1}}\right)(w) \\
& = & (\lambda\omega(e) + B^j_A\epsilon^A\wedge\omega(e)_j) (w)\\
& = & \phibl^*\Theta (w)\; .~\blob
\eeqann

\noindent Equation~(\ref{eqn:n-forms})
holds
for $n=1$ if $V_i(B,\lambda) = \lambda$ and $V_A(B,\lambda) = B_A$,
producing a parametrized version of equation~(\ref{eqn:one-form}) for
Hamiltonian particle mechanics.
Note that we cannot adapt
 Theorem~\ref{thm:nthetaV} 
to construct a connection-independent ``canonical'' $n$-form on $J^*Y$,
because
a connection is needed to define an inclusion map from $J^*Y$  into $\wedge^n Y$.

Recall from Secs.~II and IV that
$\projectable$  is in 
bijective correspondence not only
 with $T^1(Z)$ but also with
$T^1_V(\lvy)$.
The following theorem is an exact analogue of Theorem~\ref{thm:lmvecfields}.

\begin{thm}\label{thm:vecfields}
Let $v\in \projectable$. Let $\hat f_v \in \tvlvy$ correspond to
$v$ and let
$X_{\hat f_v}$ be the Hamiltonian vector field on 
$\lvy$ obtained from $\hat f_v$ via~(\ref{eqn:nkstruc}).
Let $f_v \in T^1(Z)$ also correspond to $v$ and let
$X_{f_v}$ be
the Hamiltonian vector field on $Z$ obtained from $f_v$ via~(\ref{eqn:multistruc}).
Then 
$$
\phibl_*X_{\hat f_v} =  X_{f_v}
$$
\end{thm}

\proof
From the definition of $\phi_{(B,\lambda)}$ we express its derivative
map ${\phibl}_*$ in local adapted coordinates. Let a vector field
$X$ on $\lvy$ have local adapted coordinate expression
$$
X = X^i\basisx i + X^A\basisy A + X^i_j\basispi i j + X^A_B \basispi A B
+ X^A_i\basispi A i \; .
$$
Let $\gamma(t)$ be a curve in $\lvy$ such that $\gamma(0)= w$
and $\dot\gamma(0) = X_w$.
Then
$$\phibl_*X_w = \frac{d}{dt}\phibl\circ\gamma(t) |_{t=0}\, .$$  
Using~(\ref{eqn:coordslvyz})
and the matrix identity 
$$\frac{\partial}{\partial \pi^k_l} \det \pi^r_s = \det (\pi^r_s) (\pi^{-1})^l_k$$
we can derive the local coordinate formula,
\beqa\label{eqn:hamvecphibl}
{\phibl}_*X_w
& = & X^i\basisx i + X^A\basisy A \nonumber
\\
&& \mbox{} + \det(\pi^r_s)B^j_B\left( (\pi^{-1})^i_jX^B_A - \pi^B_A(\pi^{-1})^i_k
  (\pi^{-1})^l_jX^k_l + \pi^B_A (\pi^{-1})^i_j(\pi^{-1})^l_k   X^k_l
   \right) \basisp i A \nonumber\\
&& \mbox{} + \;\det(\pi^r_s) \left(B^i_A  ( \pi^A_k (\pi^{-1})^k_j 
(\pi^{-1})^l_i X^j_l + (\pi^{-1})^k_i X^A_k
+ (\pi^{-1})^l_k (\pi^{-1})^j_i \pi^A_j X^k_l)\right. 
  \nonumber\\
  && \qquad\qquad\qquad\qquad
 + \left. \lambda(\pi^{-1})^l_k X^k_l \right)
{\basisp \, \, }.
\eeqa
Let $X = X_{\hat f_v}$ and substitute the local coordinate expressions
for the components of
$X_{\hat f_v}$ from (\ref{eqn:xflvy}) into~(\ref{eqn:hamvecphibl}).
After using the coordinate conversions~(\ref{eqn:coordslvyz}),
we 
obtain equation~(\ref{eqn:xfvcoords}), the local coordinate expression for
$X_{f_v}$ at the point $\phibl(w)$.
We must show that the resulting vector field is well-defined.
Since $di^*\theta$ and $d\hat f_v$ are tensorial, 
  the right action $R_g$ of $g \in \ga$ on $\lvy$ satisfies
$\xhat {f_v} \circ R_g = R_{g*}\xhat {f_v}$.
Now, let $w\in \lvy$ and let $g \in \ga$ satisfy
$g\cdot(B,\lambda) = (B,\lambda)$.
Then 
$\phibl(R_g(w)) = \phibl(w)$,
and consequently,
$$
\phibl_* (\xhat {f_v})_{R_g(w)} = \phibl_*(R_{g*}\xhat {f_v})_w =
\phibl_*(\xhat {f_v})_w\; .~\blob
$$

For $0\leq m \leq n+k$ define a re\-pre\-sen\-ta\-tion of $\tvlvy$ into the space of 
$\wedge^{m+1}\rnplusk$-valued
$m$-forms on $\lvy$ by 
\beq\label{eqn:badrep}
\hat f \mapsto \hat f \wedge(\iwtheta m) \, .
\eeq
The image of $\tvlvy$ under the re\-pre\-sen\-ta\-tion~(\ref{eqn:badrep}) is the set of 
 {\em degree $m$ momentum observables} on $\lvy$.
Obtain $\xhat f$ from $\hat f$ via~(\ref{eqn:nkstruc}).  By~(\ref{eqn:xflvy}),
$\xhat f \hook i^*\theta = \hat f$. 
Then by induction, 
\beq\label{eqn:268}
\xhat f \hook \iwtheta m = m\, \hat f \wedge (\iwtheta {m-1})\,  .
\eeq
By Lemma~\ref{lemma:mtheta}, for each $0 \leq m \leq n + k -2$
the form $d(\iwtheta {m+1})$ is closed and nondegenerate,
 making it a candidate for a higher degree gen\-e\-ral\-i\-za\-tion of the
$(n+k)$-symplectic form $i^*d\theta$.  However, 
using equation~(\ref{eqn:dmtheta}), we may instead use 
$i^*d\theta\wedge(\iwtheta {m})$
to eliminate a factor of  $(m+1)$ in the following generalization
of equation~(\ref{eqn:nkstruc}).
\beq\label{eqn:nkstruc-m}
d(\hat f \wedge(\iwtheta m))  = -X_{\hat f}\hook(i^*d\theta\wedge (\iwtheta {m}))\; .
\eeq
Equation~(\ref{eqn:nkstruc-m})
was derived 
using~(\ref{eqn:dmtheta}) and~(\ref{eqn:268}).

Let $\hat f, \hat g \in \tvlvy$. 
For $0 \leq m\leq n+k$ define a bracket on their images under   re\-pre\-sen\-ta\-tion~(\ref{eqn:badrep})
by
\beq\label{eqn:bracktheta}
\{\hat f \wedge(\iwtheta m ) ,\hat g \wedge(\iwtheta m )\} : =
-X_{\hat f} \hook ( X_{\hat g}\hook (i^*d\theta\wedge(\iwtheta m )) )\, .
\eeq
If $m=0$ then~(\ref{eqn:bracktheta}) defines the 
bracket on $\tvlvy$ analogous to~(\ref{eqn:Poissonn}) where $p=q=1$,
\beq\label{eqn:tvbrax}
\{\hat f, \hat g \} := X_{\hat f}(\hat g) 
=  -X_{\hat f}\hook (X_{\hat g}\hook i^*d\theta)\, .
\eeq
The vector space $\tvlvy$ is a Lie algebra under~(\ref{eqn:tvbrax}).

\begin{thm}\label{thm:fgtheta}
Let $\hat f,\hat g \in T^1_V(\lvy)$ and let $0 \leq m \leq n+k$.   Then
$$\{\hat f\wedge(\iwtheta m) ,\hat g\wedge(\iwtheta m) \}
= \{\hat f, \hat g\}\wedge(\iwtheta m) + m\,d(\hat f \wedge \hat g \wedge
(\iwtheta {m-1} ))
.$$
\end{thm}

\proof Assume $m \geq 2$. Then, using (\ref{eqn:268}), (\ref{eqn:bracktheta}) 
and (\ref{eqn:tvbrax}),
\beqann
\lefteqn{
\{\hat f\wedge(\iwtheta m) ,\hat g\wedge(\iwtheta m) \}
} \hspace{1in}\\
& = & -X_{\hat f} \hook \left(   (X_{\hat g}\hook i^*d\theta)\wedge(\iwtheta m)
     +  \,i^*d\theta\wedge (X_{\hat g}\hook  \iwtheta {m})
       \right) \\
& = & X_{\hat f}\hook \left( d\hat g \wedge (\iwtheta m))
  	+ m\,\hat g \wedge i^*(d\theta\wedge{(\wtheta {m-1})}) \right) \\
& = & X_{\hat f}(\hat g)\wedge(\iwtheta m)
      - d\hat g \wedge(X_{\hat f }\hook \iwtheta m)  \\
&& \mbox{}\quad
    + m\,\hat g \wedge (X_{\hat f}\hook i^*d\theta)\wedge (\iwtheta {m-1}) 
 	 + m\hat g\wedge i^*d\theta\wedge(X_{\hat f}\hook \iwtheta {m-1})
	\\
& = & \{\hat f ,\hat g \}\wedge(\iwtheta m)
      - m\, d\hat g\wedge \hat f \wedge(\iwtheta {m-1}) 
       - m \,\hat g\wedge d\hat f \wedge(\iwtheta {m-1}) 
\\
&& \qquad\qquad
+ m(m-1)\hat g\wedge i^*d\theta \wedge \hat f \wedge ({\iwtheta {m-2}})
\\
& = & \{\hat f , \hat g \}\wedge(\iwtheta m)
+ m\left( 
\hat f \wedge d\hat g \wedge(\iwtheta {m-1}) 
+  d\hat f \wedge \hat g \wedge {(\iwtheta {m-1})}\right.
\\
& &\mbox{} \qquad \qquad
+\left. \hat f \wedge \hat g \wedge d(\iwtheta {m-1})
\right) \\
& = & \{\hat f , \hat g \}\wedge(\iwtheta m) 
+ m\, d(\hat f\wedge \hat g \wedge (\iwtheta {m-1}))\; .
\eeqann
The proof of the $m=1$ case is a minor modification of the above argument.
The $m= 0$ case is obvious.~\blob

\bigskip

\noindent
Theorem~\ref{thm:fgtheta} is
of particular interest when $m = n-1$.  
In this case, we have reproduced for
$\lvy$ the exact analogue of the problem in equation~(\ref{eqn:pbexact}).
Just like $T^1(Z)$,
the set of 
degree $(n-1)$ momentum observables on $\lvy$
is obstructed from closing under their bracket, and the
source of the obstruction is an exact $(n-1)$-form.

Any extension of the set
of degree $(n-1)$ momentum 
observables on $\lvy$  
to include closed $\rnplusk$-valued 
$(n-1)$-forms  would have a 
representation to $HV^1(\lvy)$ with a
 kernel larger than that of the re\-pre\-sen\-ta\-tion from
$\tvlvy$ to {\em the same} vector fields.
The space $\tvlvy$ 
of observables that are both Hamiltonian and tensorial
already forms a well-defined Lie algebra under
the bracket defined in equation~(\ref{eqn:tvbrax}),
 $\tvlvy$ and $T^1(Z)$ are in bijective correspondence,
and
Theorem~\ref{thm:vecfields} dictates
  a relationship between 
the corresponding spaces of Hamiltonian  vector fields.  Thus we may 
argue  that the
$(n+k)$-symplectic geometry of $\lvy$

geometry of $Z$ for classical fields, but also possesses a Lie algebra structure for field observables that $Z$ lacks.

\section{\normalsize\bf  CONCLUSIONS}

\setcounter{equation}{0}

This investigation 
shows that the
$(n+k)$-symplectic geometry of the bundle
of vertically adapted
linear frames $\lvy$ of a field configuration space $Y$
provides a general
description of classical Hamiltonian field theories
analogous to the description of Hamiltonian par\-ti\-cle mechanics using
Norris's $n$-symplectic geometry.
A vector space of Hamiltonian
degree one tensorial
 functions  serves as the space of allowable classical field observables,
and the allowable observables  produce  Hamiltonian vector fields via the $(n+k)$-symplectic equation.

The  $(n+k)$-symplectic geometry of $\lvy$
generalizes not only the affine mul\-ti\-sym\-plec\-tic geometry
of the prototypic affine multiphase space $Z$ but also the linear mul\-ti\-sym\-plec\-tic geometry of the linear multiphase space $J^*Y$. 
This new geometry illustrates exactly why the linear
multiphase space  $J^*Y$ requires a choice  of a background
connection  in order to model classical field theories. 
Specifically,
a particular symmetry breaking of $\lvy$ is equivalent to a choice
of a connection on $Y$
which is necessary
in order to define a linear mul\-ti\-sym\-plec\-tic potential on $J^*Y$.
This symmetry breaking of $\lvy$ reinforces
the notion that the 
connection-independent affine theory on $Z$ is a more natural
starting point for mul\-ti\-sym\-plec\-tic field theories.  
We also use the  $(n+k)$-symplectic geometry of $\lvy$
to resolve a problem with
the naturally defined Poisson bracket 
of two  momentum observables on $Z$ by the construction of an 
analogous
bracket of two Hamiltonian degree one tensorial observables on $\lvy$.

In 
addition to generalizing the existing mul\-ti\-sym\-plec\-tic field theories,
our generalized symplectic geometry has an advantage over 
these theories.
In order to witness events
in a covariant field theory, we must choose
a preferred reference frame. 
The principal bundle $\lvy$ is useful because 
it describes the space of all reference frames, and sections of $\lvy$
describe preferential frames of  observers.${}^{\ref{No1}}$

The results of current work will demonstrate the utility of this generalized
symplectic geometry.
In a paper with R. O. Fulp and L. K. Norris,${}^{\ref{FLN2}}$ the author shows
that $(n+k)$-symplectic geometry on the full frame bundle $LY$ 
generalizes
the natural geometry
of {\em any} skew-symmetric tensor bundle over $Y$.
Additionally,  in a forthcoming publication,
the author will introduce momentum mappings on $\lvy$ for field theories,
setting the stage for  geometric re\-pre\-sen\-ta\-tions of frames
of field momenta and field conservation laws.

As foreshadowed by the superscript ``1'' in the notation (such as in $\tvlvy$),
the vector spaces of momentum observables may be extended to higher degree tensor fields, 
and the definitions of the Poisson bracket  may be extended to these new spaces. 
(This has been accomplished for particle mechanics, resulting in 
full Poisson and graded Poisson Schouten-Nijenhuis algebras of tensorial observables on $LM$ under the bracket in 
(\ref{eqn:Poissonn}).  See Refs.~\ref{No1} and \ref{No3}.)
Mathematically, the extended Poisson bracket would be a true Poisson bracket, acting
as a derivation on the tensor algebra. Physically,
this will create new geometrical models of Hamiltonian tensor fields,
such as  electromagnetism and Hamiltonian gravity.
Early successes in extending the algebras of observables on 
$LM$ and on $J^*Y$ (see Ref.~{\ref{Ka}})
motivate us to pursue this direction of inquiry.

Finally there may be applications of these results  
to a principal bundle formulation of
 geometric quantization of fields,
similar to that done for par\-ti\-cles.${}^{\ref{FLN1}}\,$
A prolongation of $\lvy$
may be the  appropriate setting
in which to find a faithful re\-pre\-sen\-ta\-tion of classical
observables into a space of prequantum operators, leading, for example,
to a geometric derivation of the Dirac equation.

\section*{\normalsize\bf ACKNOWLEDGEMENTS} 

 The author thanks
 R. O. Fulp and L. K. Norris for providing insights
 throughout this investigation and for reading earlier versions of this manuscript.

\bigskip\bigskip

\newcounter{references}

\begin{list}{${}^{\mbox{\small\arabic{references}}}$}{\usecounter{references}}  
\item\label{No1} L.~K. Norris, 
	``Generalized symplectic geometry on the frame bundle of a manifold,'' in
	{\em AMS Summer Research Institute on Differential Geometry} (1990),
	edited by \mbox{R.~E.} Greene and S.-T. Yau,
	{\em Proc.\ Symp.\ Pure Math.} {\bf 54}, Part 2
	(Amer.\ Math.\ Soc., Providence, RI, 1993), 435--465.
\item\label{No2} L.~K.  Norris, ``Symplectic geometry on $T^*M$ derived from
	$n$-symplectic geometry on $LM$,''  {J.~Geom.\ Phys.}\
	 {\bf 13}, 51--78 (1994). 
\item\label{Gu} C.~G\"unther, ``The polysymplectic Hamiltonian formalism in
	field theory and calculus of variations, I: The local case,''
	 {J. Diff.\ Geom.} {\bf 25}, 23--53 (1987).
\item\label{KT} J.~Kijowski and W.~M. Tulczyjew, 
	{\em A symplectic framework for field theories},
	Lect.\ Notes Phys.\ {\bf 170} (Springer-Verlag, Berlin, 1979).
\item\label{RR} R. Ragionieri and R. Ricci,
	``Hamiltonian formalism in the calculus
	of variations,''
	{Bollettino U.M.I.} {\bf 18}-B, 119--130 (1981).
\item\label{Ki2} J. Kijowski,  ``Multiphase spaces and gauge in the calculus
	of variations,'' {Bull.\ Acad.\ Sc.\ Polon.} {\bf 22}, 1219--1225 (1974).
\item\label{Ki1} J. Kijowski, ``A finite-dimensional canonical formalism in the
	classical field theory,'' {Commun.\ Math.\ Phys.} {\bf 30} , 99--128 (1973).
\item\label{Go} M. Gotay, ``A mul\-ti\-sym\-plec\-tic framework for classical field
	theory and the calculus of variations, I: Covariant Hamiltonian
	formalism,''
	in {\em Mechanics, Analysis and Geometry:  200 Years after Lagrange},
	edited by M. Francaviglia
	 (North Holland, Amsterdam, 1991), 203--235.
\item\label{GIMMsy} M. J. Gotay, J. Isenberg, J. E.  Marsden, and R.  Montgomery,   
	{\em Momentum maps and classical relativistic  fields},
	MSRI preprint (1994).
\item\label{CCI} J. F. Cari\~nena, M. Crampin, and L. A. Ibort,
	``On the multisymplectic formalism for first order field theories,''
	Diff.\ Geom.\ Appl.\ {\bf 1},  345--374 (1991).
\item\label{Ch} C. Chevalley, {\em Theory of Lie groups I},
	(Princeton University Press, Princeton, 1946).
\item\label{AM} R. Abraham and J. E. Marsden,
	{\em Foundations of Mechanics}, 2nd ed.
	(Benjamin and Cummings, New York, 1978).
\item\label{Wo}  N. Woodhouse,   
	{\em Geometric Quantization}, 2nd ed.
	(Oxford Press, Oxford  1992).
\item\label{KN} S. Kobayashi  and K.  Nomizu, 
	{\em Foundations of Differential Geometry}, Vol.\ I
	(Interscience, New York, 1963)
\item\label{Sc} J. A. Schouten, 
	``\"Uber Differentialkomitanten zweier kontravarianter Grossen,''
	Proc.\ Kon.\ Ned.\ Akad.\ Wet.\ Amsterdam {\bf 43}, 449--452 (1940).
 \item\label{Ni} A. Nijenhuis, 
	``Jacobi-type identities for bilinear differential concomitants of certain tensor
	fields,''
	{Indag.\ Math.} {\bf 17}, 390--403 (1955).
\item\label{Tr} A. Trautman, 
	``The geometry of gauge fields,''
	{Czech.\  J. Phys.\ B} {\bf 29}, 107--116 (1979).
\item\label{Bl} D. Bleecker,
	{\em Gauge Theory and Variational Principles}
	(Addison-Wesley, Reading, MA, 1981).
\item\label{Sa} G. A. Sardanashvily, {\em Gauge Theory in Jet Manifolds}
	(Hadronic Press, Palm Harbor, FL, 1993).
\item\label{FLN2} R. O. Fulp,  J. K. Lawson, L. K. Norris,  
	``Generalized symplectic geometry as a covering theory for the Hamiltonian
	theories of classical par\-ti\-cles and fields,'' 
	{J. Geom.\  Phys.} {\bf 20} 195--206 (1996).
\item\label{No3} L. K. Norris,  ``Schouten-Nijenhuis brackets,'' to appear in J.\ Math.\ Phys.
\item\label{Ka}  I. V. Kanatchikov,   ``From the Poincar\'e-Cartan form to
	a Gerstenhaber algebra of Poisson brackets in field theory,''
	in {\em Quantization, Coherent States, and Complex Structures}, 
	edited by J.-P. Antoine et al.\
	(Plenum Press, New York, 1995),  173--183.
\item\label{FLN1} R. O. Fulp, J. K.  Lawson, and  L. K. Norris,  
	``Geometric prequantization on the spin bundle based on $n$-symplectic
	geometry: The Dirac equation,''
	{Int.\ J. Theor.\ Phys.} {\bf 33}, 1011--1028 (1994).

 \end{list}

\end{document}